\documentclass{amsart}

\usepackage{amscd,amsmath,xypic,amssymb,combelow,tikz-cd,etoolbox,calligra,mathrsfs,enumitem,hyperref,mathtools}

\makeatletter
\patchcmd{\@settitle}{\uppercasenonmath\@title}{}{}{}
\makeatother

\xyoption{all}
\CompileMatrices

\makeatletter
\@addtoreset{equation}{section}
\makeatother

\newtheorem{theorem}[subsection]{Theorem}

\newtheorem{proposition}[subsection]{Proposition}
\newtheorem{lemma}[subsection]{Lemma}

\newtheorem{definition}[subsection]{Definition}

\newtheorem{remark}[subsection]{Remark}

\def\loccitt{\emph{loc. cit.}}

\def\fgl{{\mathfrak{gl}}}

\def\hgl{{\widehat{\fgl}}}

\def\fZ{{\mathfrak{Z}}}

\def\BA{{\mathbb{A}}}
\def\BC{{\mathbb{C}}}

\def\BN{{\mathbb{N}}}

\def\BP{{\mathbb{P}}}

\def\BQ{{\mathbb{Q}}}
\def\BZ{{\mathbb{Z}}}

\def\CA{{\mathcal{A}}}
\def\CB{{\mathcal{B}}}
\def\DD{{\mathcal{D}}}
\def\CE{{\mathcal{E}}}
\def\CF{{\mathcal{F}}}
\def\CG{{\mathcal{G}}}

\def\CK{{\mathcal{K}}}
\def\CL{{\mathcal{L}}}
\def\CM{{\mathcal{M}}}

\def\CO{{\mathcal{O}}}

\def\CU{{\mathcal{U}}}
\def\CV{{\mathcal{V}}}
\def\CW{{\mathcal{W}}}

\def\and{\textrm{ }\&\textrm{ }}

\def\spec{\textrm{Spec}}

\def\espec{\emph{Spec}}

\def\sym{\textrm{sym}}

\def\esym{\emph{sym}}

\def\tzeta{{\widetilde{\zeta}}}

\def\b0{{\boldsymbol{0}}}

\def\coh{\text{coh}}
\def\ecoh{\emph{coh}}

\def\kth{\text{K-th}}
\def\ekth{\emph{K-th}}

\def\ellcoh{\text{Ell}}
\def\eellcoh{\emph{Ell}}

\def\YY{Y_{t_1,t_2}(\hgl_1)}
\def\UU{U_{q_1,q_2}(\ddot{\fgl}_1)}
\def\EE{U_{q_1,q_2,p}(\ddot{\fgl}_1)}

\def\taut{\text{taut}}

\begin{document}

\title[Rational, trigonometric, elliptic algebras and sheaves on surfaces]{\Large{\textbf{Rational, trigonometric, elliptic algebras and moduli spaces of sheaves on surfaces}}}

\author[Andrei Negu\cb t]{Andrei Negu\cb t}
\address{École Polytechnique Fédérale de Lausanne (EPFL), Lausanne, Switzerland}
\address{Simion Stoilow Institute of Mathematics (IMAR), Bucharest, Romania}
\email{andrei.negut@gmail.com}

\maketitle

\begin{abstract} 

We survey the well-known Yangian of $\hgl_1$ /quantum toroidal $\fgl_1$ action on the cohomology / $K$-theory of moduli spaces of stable sheaves on surfaces, and give the generalization of this construction to elliptic cohomology. \\


\end{abstract}

$$$$

\section{Introduction}

\medskip

\subsection{} Let $H_X$ denote the singular cohomology \footnote{We employ this non-standard notation for cohomology to illustrate the fact that the discussion below applies equally well to other (co)homology theories, such as Borel-Moore homology, Chow groups, equivariant cohomology etc, with minor modifications.} of a smooth algebraic variety $X$ over $\BC$. One of the major milestones in geometric representation theory was the construction (independently due to Grojnowski \cite{grojnowski} and Nakajima \cite{nakajima}) of an action 
\begin{equation}
\label{eqn:hilbert}
\hgl_1 \curvearrowright \bigoplus_{n=0}^{\infty} H_{S^{[n]}}
\end{equation}
In the left-hand side, we have the infinite-dimensional Heisenberg Lie algebra. In the right-hand side, the Hilbert schemes $S^{[n]}$ of $n$ points on a smooth projective surface $S$ are perhaps some of the simplest moduli spaces of coherent sheaves. They parameterize length $n$ subschemes of $S$, or equivalently, colength $n$ ideal sheaves. The latter interpretation allows one to think of Hilbert schemes as moduli spaces of rank 1 coherent sheaves on $S$, and therefore suggests that \eqref{eqn:hilbert} should admit a higher rank generalization. Indeed, one may consider the moduli spaces
\begin{equation}
\label{eqn:m intro}
\CM = \bigsqcup_{c_2 \in \BZ} \CM_{(r,c_1,c_2)}
\end{equation}
parameterizing stable sheaves $\CF$ on $S$ of fixed rank $r$ and first Chern class $c_1$, but variable second Chern class $c_2$ (see Subsection \ref{sub:basic moduli} for an overview of $\CM$). In this direction, Baranovsky showed in \cite{baranovsky} that one can generalize \eqref{eqn:hilbert} to an action 
\begin{equation}
\label{eqn:moduli}
\hgl_1 \curvearrowright H_{\CM} = \bigoplus_{c_2 \in \BZ} H_{\CM_{(r,c_1,c_2)}}
\end{equation}
Another direction in which one may seek to generalize \eqref{eqn:hilbert} is to enlarge the infinite-dimensional Heisenberg Lie algebra. Inspired by work of \cite{gv, nak quiver, varagnolo} on quiver varieties, the corresponding enlargement is an algebra known as the Yangian
$$
\YY \supset \hgl_1
$$
(see \cite{drinfeld} for Yangians of finite type Lie algebras, of which the object above is a natural generalization) defined over the ring $\BZ[t_1,t_2]$. We will recall the definition of the Yangian in Subsection \ref{sub:yang} and prove the following generalization of \eqref{eqn:moduli}.

\medskip

\begin{theorem}
\label{thm:coh intro}

(Theorem \ref{thm:coh}) For any smooth projective surface $S$, and $r,c_1$ satisfying Assumptions A and S in \eqref{eqn:assumption a} and \eqref{eqn:assumption s} (respectively), there is an action
\begin{equation}
\label{eqn:moduli yangian}
\YY \curvearrowright H_{\CM}
\end{equation}

\end{theorem}

\medskip

\noindent Above, one must be careful to properly define the notion of an ``action" of the Yangian on the cohomology of moduli spaces of stable sheaves. In Definition \ref{def:action}, the appropriate notion will be shown to be that of an abelian group homomorphism
\begin{equation}
\label{eqn:moduli yangian explicit}
\YY \rightarrow \text{Hom} \left( H_{\CM}, H_{\CM \times S} \right)
\end{equation}
satisfying several axioms; the parameters $t_1,t_2$ are identified with the Chern roots of the cotangent bundle of $S$, in a sense that will be made precise in Definition \ref{def:action}.

\medskip

\noindent If $S = \BA^2$ and one works with the moduli space of framed sheaves (a close relative of $\CM$, which we will recall in Subsection \ref{sub:framed}) then the cohomology of $S$ is trivial. In this case, \eqref{eqn:moduli yangian explicit} yields an honest action of the Yangian on the cohomology groups of the moduli spaces of framed sheaves. Since the latter spaces are isomorphic to Nakajima quiver varieties for the Jordan quiver, the $S=\BA^2$ analogue of \eqref{eqn:moduli yangian} (see Theorem \ref{thm:yang}) can be construed as the Jordan quiver version of \cite{gv, nak quiver, varagnolo}.

\medskip

\subsection{}

In representation theory, there are three ways to ``affinize" Lie algebras
\begin{equation}
\label{eqn:affinization}
\Big( \text{Yangians} \Big) \leftarrow \Big( \text{quantum loop groups} \Big) \leftarrow \Big( \text{elliptic quantum groups} \Big)
\end{equation}
that correspond to solving the Yang-Baxter equation with parameter in $\BC$, $\BC^*$ and an elliptic curve $E$, respectively (thus the algebras above are usually called rational, trigonometric and elliptic, respectively). The arrows in \eqref{eqn:affinization} mean that by suitably degenerating elliptic quantum groups, one may obtain quantum loop groups, and by suitably degenerating quantum loop groups, one may obtain Yangians. 

\medskip

\noindent On the geometric side, we also have three important oriented homology theories
\begin{equation}
\label{eqn:homology}
\Big( \text{cohomology} \Big) \leftarrow \Big( K\text{-theory} \Big) \leftarrow  \Big( \text{elliptic cohomology} \Big)
\end{equation}
As we have already seen in the previous Subsection, the Yangian of $\hgl_1$ acts on the cohomology of moduli spaces of stable sheaves. Therefore, the following is a natural analogue of Theorem \ref{thm:coh intro}. Below and henceforth, we let $K_X$ denote either the 0-th topological or algebraic $K$-theory group of a smooth projective variety $X$ over $\BC$.

\medskip

\begin{theorem}
\label{thm:k-theory intro}

(Theorem \ref{thm:k-theory}) For any smooth projective surface $S$, and $r,c_1$ satisfying Assumptions A and S in \eqref{eqn:assumption a} and \eqref{eqn:assumption s} (respectively), there is an action
\begin{equation}
\label{eqn:moduli toroidal}
\UU \curvearrowright K_{\CM}
\end{equation}

\end{theorem}

\medskip

\noindent We will recall the definition of the algebra $\UU$ in Subsection \ref{sub:tor}. It is called quantum toroidal $\fgl_1$, as well as the Ding-Iohara-Miki algebra (\cite{di,miki}), and it was studied in numerous works (closest to our point of view being \cite{fhhsy, fjmm}). As an algebra, it is defined over the ring $\BZ[q_1^{\pm 1}, q_2^{\pm 1}]$, where the formal parameters $q_1,q_2$ should be construed as the exponentials of the parameters $t_1,t_2$ of the previous Subsection. From a geometric point of view, $q_1+q_2$ and $q_1q_2$ are identified with the $K$-theory classes of the cotangent bundle and the canonical line bundle of $S$, respectively.

\medskip

\noindent Although it is not apparent from the usual definition (which we will recall in Subsection \ref{sub:tor}), there is an inclusion of algebras
\begin{equation}
\label{eqn:q heisenberg}
\Big(q\text{-Heisenberg algebra}\Big) \hookrightarrow \UU
\end{equation}
Because of this, \eqref{eqn:moduli toroidal} yields an action of the $q$-Heisenberg algebra on the $K$-theory of moduli spaces of stable sheaves, i.e. a $q$-version of the action \eqref{eqn:moduli} \footnote{However, the tools surveyed in the present paper do not give a geometric incarnation of the action \eqref{eqn:q heisenberg}, a task for which one should appeal to the correspondences developed in \cite{N mod, N hecke}. Doing so naturally leads to an action of the so-called elliptic Hall algebra from \cite{bs}, a close relative of quantum toroidal $\fgl_1$, on the $K$-theory of moduli spaces of stable sheaves.}. 

\medskip

\subsection{} The novelty in the present paper is to run the above program in the setting of elliptic cohomology. We will follow the axiomatic viewpoint of \cite{gkv}, which associates to an elliptic curve $E$ the contravariant functor of (0-th) elliptic cohomology
\begin{equation}
\label{eqn:elliptic intro}
\Big(\text{smooth varieties over }\BC \Big) \quad \xrightarrow{X \mapsto \ellcoh_X} \quad \text{Commutative rings}
\end{equation}
and regards Chern classes of rank $r$ vector bundles $\CV$ on $X$ as maps of schemes
\begin{equation}
\label{eqn:chern intro}
\spec(\ellcoh_X) \xrightarrow{c_{\CV}} E^{(r)}
\end{equation}
where $E^{(r)} = E^r/S_r$ denotes the $r$-th symmetric power of $E$. Naively, one would imagine elliptic cohomology classes to be elements of $\ellcoh_X$, i.e. regular functions on $\spec(\ellcoh_X)$. However, the correct thing to do would be to broaden one's understanding and declare that elliptic cohomology classes on $X$ are
\begin{multline*}
\text{sections of line bundles on }\spec(\ellcoh_X) \Leftrightarrow \\ \Leftrightarrow \text{elements of locally free rank 1 modules over }\ellcoh_X
\end{multline*}
Indeed, a big source of elliptic cohomology classes are the pull-backs of various meromorphic functions on $E^{(r)}$ to $\spec(\ellcoh_X)$ via the map $c_{\CV}$ of \eqref{eqn:chern intro}. As these meromorphic functions are naturally interpreted as sections of line bundles on $E^{(r)}$, their pull-backs via $c_{\CV}$ are also naturally sections of line bundles on $\spec(\ellcoh_X)$. 

\medskip

\subsection{} For any map of smooth varieties $f : X \rightarrow Y$, one has a ring homomorphism
$$
\ellcoh_Y \xrightarrow{f^*} \ellcoh_X
$$
For any proper map $f : X \rightarrow Y$, push-forward is defined by \cite{gkv} as a map 
$$
\Big\{ \text{a locally free rank 1 }\ellcoh_X\text{-module} \Big\} \xrightarrow{f_*} \ellcoh_Y
$$
of $\ellcoh_Y$-modules, for any proper morphism $f : X \rightarrow Y$. However, as we will recall in Subsection \ref{sub:push 1}, the above locally free rank 1 module is completely determined by the virtual relative tangent bundle of the morphism $f$. Moreover, in the present paper we will only encounter two situations: when $f$ is a regular embedding and when $f$ is a projective bundle. In both of these cases, there exist explicit formulas for $f_*$ which completely encode the locally free rank 1 module in question. Thus, we will often abuse notation and write
\begin{equation}
\label{eqn:wrong}
f_* : \ellcoh_X \rightarrow \ellcoh_Y
\end{equation}
throughout the present paper. The reader may then easily deduce which locally free rank 1 $\ellcoh_X$-module needs to go in the domain of $f_*$ just by looking at the formula for $f_*$ (see Subsection \ref{sub:push 3} for a defining example of this principle). Thus, we are not really sacrificing mathematical precision in writing formulas such as \eqref{eqn:wrong}; we are simply sacrificing notational precision, and the benefit is increased legibility.

\medskip

\subsection{}

We will consider the elliptic cohomology of the moduli space \eqref{eqn:m intro} of stable sheaves, for any fixed $r \in \BN$ and $c_1 \in H^2(S,\BZ)$ that satisfy Assumptions A and S. In formulas \eqref{eqn:e elliptic series}, \eqref{eqn:f elliptic series}, \eqref{eqn:h elliptic series}, we construct formal series of operators
\begin{equation}
\label{eqn:series intro}
e(z), f(z), h^\pm(z) : \ellcoh_{\CM} \rightarrow \ellcoh_{\CM \times S}[[z,z^{-1}]]
\end{equation}
Above, the symbol $z$ plays the role of the usual variable on the cover
$$
\BC^* \rightarrow \BC^*/p^{\BZ} = E
$$
where $p$ is a complex number with absolute value strictly contained between 0 and 1. Our main result is the commutation relations between the series of operators \eqref{eqn:series intro}. These will be presented in Theorem \ref{thm:elliptic}, and they correspond to an action
\begin{equation}
\label{eqn:intro elliptic}
\EE \ ``\curvearrowright" \ \ellcoh_{\CM}
\end{equation}
where the object in the left-hand side is the elliptic quantum toroidal algebra defined in \cite{ko} (with central charge equal to 1 and without the cubic relation of \loccitt, see Subsection \ref{sub:ko}). The reason we do not state \eqref{eqn:intro elliptic} as a Theorem is of a rather pedantic nature: the left-hand side is an algebra defined over power series in $p$, while the right-hand side is defined for a specific value of $p \in \BC^*$. Thus, we leave it as an open (and very interesting) question to adapt the contents of the present paper from the setting of fixed $p$ to that of variable $p$, especially by working with elliptic cohomology over the moduli space of oriented elliptic curves as in 
\cite{lurie}.

\medskip 

\subsection{} I would like to thank the organizers of the 2023 MSJ-SI ``Elliptic Integrable Systems, Representation Theory and Hypergeometric Functions" workshop at the University of Tokyo (with special gratitude to Hitoshi Konno) for the invitation and 
inspiration which led to the writing of the present paper. I gratefully acknowledge support from the NSF grant DMS-1845034, the MIT Research Support Committee and the PNRR grant CF 44/14.11.2022 titled ``Cohomological Hall algebras of smooth surfaces and applications".

\medskip

\section{Preliminaries on moduli of sheaves on surfaces}
\label{sec:moduli}

\medskip

\subsection{} 
\label{sub:basic moduli}

Fix a smooth projective connected surface $S$ over an algebraically closed field of characteristic 0 (we will henceforth abuse notation and denote the field by $\BC$) and a ample divisor $H \subset S$. A coherent sheaf $\CF$ on $S$ is called stable if all proper subsheaves $\CG \subsetneq \CF$ have reduced Poincar\'e polynomial (defined with respect to $H$) smaller than that of $\CF$; see \cite[Section 5]{N shuffle surf} and the references therein for an overview. 

\medskip

\begin{theorem} \label{thm:stable} For any $(r,c_1,c_2) \in \BN \times H^2(S,\BZ) \times \BZ$, there exists a moduli space $\CM_{(r,c_1,c_2)}$ corepresenting the functor of stable sheaves $\CF$ on $S$ of rank $r$, first Chern class $c_1$ and second Chern class $c_2$.

\end{theorem}

\medskip

\noindent Bogomolov's inequality states that
\begin{equation}
\label{eqn:bogomolov}
\CM_{(r,c_1,c_2)} = \varnothing \quad \text{if} \quad c_2 < \frac {r-1}{2r} c_1^2
\end{equation}
Throughout the present paper, we fix $r,c_1$ and write
\begin{equation}
\label{eqn:disjoint union}
\CM = \bigsqcup_{c_2 \in \BZ} \CM_{(r,c_1,c_2)}
\end{equation}
although, as mentioned in \eqref{eqn:bogomolov}, the number $c_2$ is actually bounded below. We will henceforth make the following assumptions, as in \cite{N shuffle surf}. The first of these states that
\begin{equation}
\label{eqn:assumption a}
\textbf{Assumption A: } \gcd(r,c_1\cdot H) = 1
\end{equation}
and implies that $\CM$ is projective (because Assumption A implies that any seminstable sheaf is stable), and that it also represents the functor referenced in Theorem \ref{thm:stable}. In other words, this means that there exists a universal sheaf
\begin{equation}
\label{eqn:universal}
\xymatrix{ \CU \ar@{.>}[d] \\ \CM \times S}
\end{equation}
whose restriction to any point $\CF \in \CM$ is isomorphic to $\CF$ as a sheaf on $S$. The second assumption pertains to the canonical bundle $\CK_S$ of $S$, and it reads
\begin{equation}
\label{eqn:assumption s}
\textbf{Assumption S:} \text{ either } \CK_S \cong \CO_S \text { or } c_1(\CK_S) \cdot H < 0 
\end{equation}
The latter assumption implies that $\CM$ is smooth, and although it will be in force throughout the present paper, it can be sidestepped if one is prepared to work with derived moduli spaces of sheaves. However, Assumption A is quite crucial for us, because we will fundamentally use the universal sheaf $\CU$.

\medskip

\subsection{}

If $\CF, \CF'$ are sheaves on $S$ such that there exists a short exact sequence
\begin{equation}
\label{eqn:ses}
0 \rightarrow \CF' \rightarrow \CF \rightarrow \BC_x \rightarrow 0
\end{equation}
for some closed point $x \in S$ (where $\BC_x$ denotes the skyscraper sheaf at $x$), then we will say that $\CF$ and $\CF'$ are Hecke modifications of each other and write this as
$$
\CF' \subset_x \CF
$$
An important consequence of Assumption A is the fact (see \cite[Proposition 5.5]{N shuffle surf}) that $\CF$ is stable if and only if $\CF'$ is stable. Thus, we may define the moduli space
\begin{equation}
\label{eqn:simple}
\fZ = \Big\{(\CF,\CF') \text{ both stable, such that } \CF' \subset_x \CF \Big\}
\end{equation}
which is smooth and projective (\cite[Proposition 2.10]{N shuffle surf}). Let us consider the maps
\begin{equation}
\label{eqn:simple diag}
\xymatrix{ & \fZ \ar[dl]_{\pi_+} \ar[d]^{\pi_S} \ar[dr]^{\pi_-} & \\
\CM' & S & \CM} 
\end{equation}
which remember $\CF',x,\CF$ (respectively) in the notation of \eqref{eqn:ses}. Above, $\CM$ and $\CM'$ are simply two copies of the moduli space \eqref{eqn:disjoint union}, although one should think that the grading by $c_2$ in the moduli space $\CM'$ is always 1 more than that of $\CM$. 

\medskip

\noindent If $\CV$ is a locally free sheaf on a variety $X$, then we define its projectivization as 
$$
\BP_X(\CV) = \text{Proj}_X(\text{Sym}^\bullet \CV)
$$
More generally, if $\CE$ is a coherent sheaf of projective dimension 1, its projectivization is defined in \eqref{eqn:proj sheaf}. With this in mind, we have the following description of $\fZ$.

\medskip

\begin{proposition}
\label{prop:projective bundles}

(\cite[Propositions 2.8 and 2.10]{N shuffle surf}) The maps $\pi_- \times \pi_S$ and $\pi_+ \times \pi_S$ realize $\fZ$ as the projectivizations
\begin{equation}
\label{eqn:projectivizations}
\BP_{\CM \times S}(\CU) \quad \text{and} \quad \BP_{\CM' \times S}({\CU'}^\vee[1] \otimes \CK_S)
\end{equation}
where $\CU$ and $\CU'$ denote the universal sheaves on $\CM \times S$ and $\CM' \times S$, respectively.

\end{proposition}

\medskip

\noindent Consider the line bundle 
\begin{equation}
\label{eqn:tautological}
\xymatrix{ \CL \ar@{.>}[d] \\ \fZ}
\end{equation}
whose fiber over a point $(\CF', \CF)$ is the one-dimensional vector space $\Gamma(S,\CF/\CF')$. Under the isomorphisms \eqref{eqn:projectivizations}, we have $\CL \cong \CO(1)$ and $\CL \cong \CO(-1)$, respectively, where $\CO(1)$ denotes the tautological quotient line bundle on any projectivization. In geometric terms, if a point of $\BP_X(\CV)$ parameterizes a one-dimensional quotient vector space $\CV_x \twoheadrightarrow Q$ for some $x \in X$, then the fiber of $\CO(1)$ over such a point is the one-dimensional vector space $Q$ itself.

\medskip

\subsection{}
\label{sub:framed} 

When $S = \BA^2$, the moduli space of stable sheaves is not well-defined. However, there is an alternate provided by the so-called moduli space of framed sheaves on $\BP^2$. Intuitively, this is because framing (like stability) gives a way to control the automorphisms of sheaves, leading to well-behaved moduli spaces.

\medskip

\begin{definition}
\label{def:moduli}

Consider $\BP^2 = \BA^2 \sqcup \infty$ and define the moduli space
\begin{multline}
\label{eqn:framed}
\CM_{\BA^2} = \Big \{\text{rank } r \text{ torsion-free sheaves }\CF \text{ on }\BP^2, \text{ locally free near }\infty, \\ \text{together with an isomorphism } \CF|_{\infty} \cong \CO_{\infty}^{\oplus r} \Big\}
\end{multline}

\medskip

\end{definition}

\noindent The moduli space $\CM_{\BA^2}$ has connected components indexed by $n = c_2(\CF)$, namely
$$
\CM_{\BA^2} = \bigsqcup_{n = 0}^{\infty} \CM_{\BA^2,n}
$$
For any $n$, the connected component $\CM_{\BA^2,n}$ is a smooth quasiprojective algebraic variety, and can be explicitly presented (via the famous Atiyah-Drinfeld-Hitchin-Manin construction) as the space of quadruples of linear maps
\begin{multline}
\label{eqn:adhm}
\Big\{(X,Y,A,B) \text{ s.t. } \BC^n \xrightleftharpoons[Y]{X} \BC^n, \BC^r \xrightleftharpoons[B]{A} \BC^n, [X,Y] + AB = 0 \Big\}^{\text{stable}} \Big / GL_n
\end{multline}
Above, a quadruple $(X,Y,A,B)$ is called stable if there is no proper subspace of $\BC^n$ which contains $\text{Im }A$ and is preserved by both $X$ and $Y$. The action of $GL_n$ on quadruples is given by conjugating $X,Y$ and multiplying $A,B$ on either side. 

\medskip

\begin{remark} 

Although $\CM_{\BA^2}$ parameterizes sheaves $\CF$ on $\BP^2$, such sheaves are trivialized on the divisor $\infty$ and so all their interesting behavior happens on $\BA^2$. Thus, we will often restrict such $\CF$ to a sheaf on $\BA^2$ in what follows. In particular, the analogue of the universal sheaf \eqref{eqn:universal} is defined as
$$
\xymatrix{ \CU \ar@{.>}[d] \\ \CM_{\BA^2} \times \BA^2}
$$
Similarly, the natural analogue of the moduli space \eqref{eqn:simple} requires the point $x$ to be away from $\infty$, so we put $S = \BA^2$ in diagram \eqref{eqn:simple diag}. We will also encounter the restriction of the universal sheaf to $\CM_{\BA^2} \times \{\text{origin}\}$
\begin{equation}
\label{eqn:restricted universal}
\xymatrix{ \CU_{\circ} \ar@{.>}[d] \\ \CM_{\BA^2}}
\end{equation}
\end{remark}

\medskip

\noindent The dilation torus action $T = \BC^* \times \BC^* \curvearrowright \BA^2$ induces an action
$$
T \curvearrowright \CM_{\BA^2}
$$
In the language of \eqref{eqn:adhm}, the torus $T$ scales the matrices $X,Y,A,B$ by the characters $\frac 1{q_1},\frac 1{q_2},1,\frac 1{q_1q_2}$, respectively, where $q_1,q_2$ are the elementary characters of $T$. 

\bigskip

\section{Cohomology}
\label{sec:coh}

\medskip

\subsection{} 
\label{sub:yang} 

We start with a warm-up on the Yangian of $\hgl_1$ and its action on the equivariant cohomology of moduli spaces of framed sheaves on $\BA^2$. We will work over the ring of polynomials in two variables $t_1$ and $t_2$. Let us write $t = t_1+t_2$ and consider the rational functions 
\begin{equation}
\label{eqn:def zeta yang}
\zeta^{\BC}(x) = \frac {(x+t_1)(x+t_2)}{x(x+t)}
\end{equation}
\begin{equation}
\label{eqn:def zeta yang tilde}
\tzeta^{\BC}(x) = \zeta^{\BC}(x)(x+t)(x-t) = \frac {(x+t_1)(x+t_2)(x-t)}{x}
\end{equation}
The following definition is inspired by the Yangians of finite type Lie algebras that were constructed in \cite{drinfeld}, and studied in numerous works afterwards.

\medskip

\begin{definition}
\label{def:yang}

(\cite{tsymbaliuk}) The Yangian of $\hgl_1$ is the algebra
$$
\YY = \BZ[t_1,t_2] \Big\langle e_n, f_n, h_n \Big \rangle_{n \geq 0} \Big / \text{relations \eqref{eqn:rel 1 yang}--\eqref{eqn:rel 5 yang}}
$$
The defining relations are best written in terms of the generating series
$$
e(z) = \sum_{n=0}^{\infty} \frac {e_n}{z^{n+1}}, \qquad f(z) = \sum_{n=0}^{\infty} \frac {f_n}{z^{n+1}}, \qquad h(z) = 1 + \sum_{n=0}^{\infty} \frac {h_n}{z^{n+1}}
$$
and take the form
\begin{equation}
\label{eqn:rel 1 yang}
\left[ e(z)e(w) \tzeta^{\BC}(w-z) \right]_{z^{<0},w^{<0}} = \left[ e(w)e(z)  \tzeta^{\BC}(z-w) \right]_{z^{<0},w^{<0}}
\end{equation}
\begin{equation}
\label{eqn:rel 2 yang}
\left[ f(w)f(z)  \tzeta^{\BC}(w-z) \right]_{z^{<0},w^{<0}} = \left[ f(z)f(w) \tzeta^{\BC}(z-w) \right]_{z^{<0},w^{<0}}
\end{equation}
\begin{equation}
\label{eqn:rel 3 yang}
h(z)e(w) =  \left[ e(w)h(z) \frac {\zeta^{\BC}(z-w)}{\zeta^{\BC}(w-z)} \right]_{z \gg w|z^{\leq 0},w^{<0}}
\end{equation}
\begin{equation}
\label{eqn:rel 4 yang}
f(w) h(z) = \left[ h(z) f(w) \frac {\zeta^{\BC}(z-w)}{\zeta^{\BC}(w-z)} \right]_{z \gg w|z^{\leq 0},w^{<0}}
\end{equation}
\begin{equation}
\label{eqn:rel 5 yang}
\left [f(z), e(w) \right] = \frac {t_1t_2}t \cdot \frac {h(z)-h(w)}{z-w}
\end{equation}
as well as $[h(z),h(w)] = 0$. 

\end{definition}

\medskip

\noindent Let us explain the meaning of the square brackets in \eqref{eqn:rel 1 yang}--\eqref{eqn:rel 5 yang}, and how to turn them into equalities of symbols: in formulas \eqref{eqn:rel 1 yang} and \eqref{eqn:rel 2 yang}, we cancel out the factor $z-w$ from the denominator and equate the coefficients of all $\{z^aw^b\}_{a<0,b<0}$ in the left and right-hand sides. In formulas \eqref{eqn:rel 3 yang} and \eqref{eqn:rel 4 yang}, we expand the rational function in the RHS in non-negative powers of $w/z$, and equate the coefficients of all $\{z^aw^b\}_{a\leq 0, b <0}$ in the left and right-hand sides. Finally, in formula \eqref{eqn:rel 5 yang}, we simply equate the coefficients of all $\{z^aw^b\}_{a<0,b<0}$ in the left and right-hand sides.

\medskip

\begin{remark}
\label{rem:extra relation yang}

In many works (such as \cite{tsymbaliuk}), the algebra $\YY$ involves the additional cubic relations analogous to the Drinfeld-Serre relations in finite types
\begin{equation}
\label{eqn:cubic yang}
\sum_{\sigma \in S_3} [e_{n_{\sigma(1)}},[e_{n_{\sigma(2)}},e_{n_{\sigma(3)}+1}]]  =  \sum_{\sigma \in S_3} [f_{n_{\sigma(1)}},[f_{n_{\sigma(2)}},f_{n_{\sigma(3)}+1}]]  =  0
\end{equation}
for all $n_1,n_2,n_3 \geq 0$. These additional relations hold in all the modules of geometric nature considered in the present paper, namely the cohomology groups of moduli spaces of stable/framed sheaves (this was proved at the level of $K$-theory in \cite{N hecke}, and the version in cohomology can be deduced via the Chern character isomorphism). 

\end{remark}

\medskip

\subsection{}

Let us unpack relations \eqref{eqn:rel 1 yang}--\eqref{eqn:rel 5 yang}. The last of these relations is the easiest one, as it yields the following identity for all $m,n \geq 0$
\begin{equation}
\label{eqn:rel 5 yang explicit}
[f_n,e_m] = - \frac {t_1t_2}t \cdot h_{n+m}
\end{equation}
To express relations \eqref{eqn:rel 3 yang} and \eqref{eqn:rel 4 yang}, we need to consider the power series expansion
\begin{equation}
\label{eqn:expansion yangian}
\frac {\zeta^{\BC}(z-w)}{\zeta^{\BC}(w-z)} = \frac {(z-w+t_1)(z-w+t_2)(z-w-t)}{(z-w-t_1)(z-w-t_2)(z-w+t)} = 1 + \sum_{a=3}^{\infty} \sum_{b=0}^{a-2} t_1t_2 \gamma_{ab} \frac {w^b}{z^a}
\end{equation}
for various polynomial expressions $\gamma_{ab}$ in $t=t_1+t_2$ and $t_1t_2$. The fact that we can extract a factor of $t_1t_2$ from the coefficients in the right-hand side of \eqref{eqn:expansion yangian} will be very important for our geometric constructions in the subsequent Subsections, and it is due to the fact that the left-hand side of \eqref{eqn:expansion yangian} is equal to 1 when either $t_1=0$ or when $t_2=0$. With this in mind, relations \eqref{eqn:rel 3 yang}--\eqref{eqn:rel 4 yang} read
\begin{align}
&[h_n,e_m] = t_1t_2 \sum_{a=3}^{\infty} \sum_{b=0}^{a-2} \gamma_{ab} \cdot e_{m+b}h_{n-a} \label{eqn:rel 3 yang explicit} \\
&[f_m,h_n] = t_1t_2 \sum_{a=3}^{\infty} \sum_{b=0}^{a-2} \gamma_{ab} \cdot h_{n-a} f_{m+b} \label{eqn:rel 4 yang explicit}
\end{align}
respectively, for all $m,n \geq 0$ (in the right-hand sides of the expressions above, we make the convention that $h_{-1}=1$ and $h_{-2} = h_{-3} = \dots = 0$).

\medskip

\noindent Finally, to explicitly write down relations \eqref{eqn:rel 1 yang}--\eqref{eqn:rel 2 yang}, we observe that
$$
\widetilde{\zeta}^{\BC}(z-w) = (z-w)^2 - (t_1^2 + t_1t_2 + t_2^2) - \frac {t_1t_2t}{z-w}
$$
Therefore, relations \eqref{eqn:rel 1 yang}--\eqref{eqn:rel 2 yang} respectively imply the following for all $m,n \geq 0$
\begin{multline}
\label{eqn:rel 1 yang explicit}
[e_{n+3},e_m] - 3 [e_{n+2},e_{m+1}] + 3 [e_{n+1},e_{m+2}] - [e_n,e_{m+3}] - \\ - (t_1^2+t_1t_2+t_2^2)([e_{n+1},e_m]-[e_n,e_{m+1}]) + t_1t_2t(e_ne_m+e_me_n) = 0
\end{multline}
\begin{multline}
\label{eqn:rel 2 yang explicit}
[f_m,f_{n+3}] - 3 [f_{m+1},f_{n+2}] + 3 [f_{m+2},f_{n+1}] - [f_{m+3},f_n] - \\ - (t_1^2+t_1t_2+t_2^2)([f_{m},f_{n+1}]-[f_{m+1},f_n]) + t_1t_2t(f_mf_n+f_nf_m) = 0
\end{multline}

\medskip

\begin{remark}
\label{rem:quotients}

To make the connection with geometry completely rigorous, we need to make two small modifications to Definition \ref{def:yang}. First of all, since all relations \eqref{eqn:rel 1 yang}--\eqref{eqn:rel 5 yang} are symmetric in $t_1$ and $t_2$, we may define 
$$
\YY \text{ as a }\BZ[t_1,t_2]^{\emph{sym}}\text{-algebra}
$$
Secondly, we need to ensure that 
\begin{equation}
\label{eqn:divisible}
[x,y] \ \text{is divisible by} \ t_1t_2
\end{equation}
for all $x,y \in \YY$. Relations \eqref{eqn:rel 5 yang explicit}, \eqref{eqn:rel 3 yang explicit} and \eqref{eqn:rel 4 yang explicit} already imply that any commutator $[f_n,e_m]$, $[h_n,e_m]$, $[f_m,h_n]$ is a multiple of $t_1t_2$, but relations \eqref{eqn:rel 1 yang explicit}, \eqref{eqn:rel 2 yang explicit} are unfortunately not strong enough to imply that the commutators $[e_n,e_m]$ and $[f_n,f_m]$ are multiples of $t_1t_2$. Therefore, we must explicitly add the symbols
\begin{equation}
\label{eqn:formally adding}
\frac {[e_{n_1},\dots,[e_{n_{k-1}},[e_{n_k},e_{n_{k+1}}]] \dots]}{(t_1t_2)^k} \quad \text{and} \quad \frac {[f_{n_1},\dots,[f_{n_{k-1}},[f_{n_k},f_{n_{k+1}}]] \dots]}{(t_1t_2)^k}
\end{equation}
to the algebra $\YY$, for all $n_1,\dots,n_{k+1} \geq 0$. Naturally, the symbols \eqref{eqn:formally adding} must satisfy the Leibniz rule and Jacobi identities in the natural sense when multiplied or commuted with each other or with the generators $e_n,f_n,h_n$ (\cite[Section 4.10]{N hecke}).

\end{remark}

\medskip

\subsection{}
\label{sub:cohomology}

For a smooth algebraic variety $X$ over $\BC$, its cohomology
\begin{equation}
\label{eqn:singular cohomology}
H_X = H^*(X,\BZ)
\end{equation}
will refer to the singular cohomology ring with integer coefficients. However, all results in the present Section also hold for other cohomology theories, such as Borel-Moore homology or Chow rings. In fact, the only properties of cohomology that we will need is the existence of pull-back maps for local complete intersection morphisms, push-forward maps for proper morphisms, and Chern classes
\begin{equation}
\label{eqn:chern}
c_i(\CV) \in H_X
\end{equation}
for any locally free sheaf $\CV$ on $X$. These may be assembled in the Chern polynomial
$$
c(\CV,z) = \sum_{i=0}^{v} z^{v-i} (-1)^i c_i(\CV) \in H_X[z]
$$
where $v = \text{rank }\CV$. The Chern polynomial only depends on the $K$-theory class of $\CV$. This means that if we have a short exact sequence of coherent sheaves
\begin{equation}
\label{eqn:ses general}
0 \rightarrow \CW \rightarrow \CV \rightarrow \CE \rightarrow 0
\end{equation}
with $\CV,\CW$ locally free of ranks $v,w$, then we may define
$$
c(\CE,z) = \frac {c(\CV,z)}{c(\CW,z)}
$$
The expression above does not depend on the choice of short exact sequence \eqref{eqn:ses general}, so it allows one to unambiguosly define the Chern classes of $\CE$ itself by
$$
c(\CE,z) = \sum_{i=0}^{\infty} z^{v-w-i} (-1)^i c_i(\CE) \in H_X ((z^{-1}))
$$

\medskip

\noindent We will also encounter equivariant cohomology in the present paper. For a smooth algebraic variety $X$ endowed with an action of a torus $T$, the ring
$$
H_T^*(X)
$$
enjoys the same functoriality properties as usual cohomology \eqref{eqn:singular cohomology}. The two main differences between these two rings is that equivariant cohomology is an algebra over the polynomial ring
$$
H_T^*(\text{point}) = \BZ[\text{Lie}(T)]
$$
and that the grading $*$ is not bounded above by the real dimension of $X$ anymore (for instance, the grading on $H_T^*(\text{point})$ has linear functions on $\text{Lie}(T)$ in degree 2).

\medskip

\subsection{}
\label{sub:correspondences}

All the geometric operators that we will define in the present paper are given by correspondences. In the most general form of this concept, a correspondence between smooth projective algebraic varieties $X$ and $Y$ is a 
class
$$
\Gamma \in H_{X \times Y}
$$
which induces the operator
$$
\Phi_\Gamma : H_Y \xrightarrow{p_Y^*} H_{X \times Y} \xrightarrow{\smile \Gamma} H_{X \times Y} \xrightarrow{p_{X*}} H_X 
$$
where $p_X : X \times Y \rightarrow X$ and $p_Y : X \times Y \rightarrow Y$ denote the standard projections. There is a natural notion of composition of correspondences (whose exact definition we leave as an exercise to the reader)
$$
(\Gamma,\Gamma') \in H_{X \times Y} \times H_{Y \times Z} \quad \leadsto \quad \Gamma \circ \Gamma' \in H_{X \times Z}
$$
such that
$$
\Phi_\Gamma \circ \Phi_{\Gamma'} = \Phi_{\Gamma \circ \Gamma'}
$$
With this in mind, we observe that while all the results in the present paper are stated as equalities of compositions of the various operators $\Phi_\Gamma$, this is only for convenience. In fact, all our results hold as equalities of compositions of correspondences $\Gamma$. This is important from a technical point of view, since the assignment
\begin{equation}
\label{eqn:assingment correspondences}
\Gamma \leadsto \Phi_\Gamma
\end{equation}
is only injective for cohomology theories which satisfy the K\"unneth decomposition
$$
H_{X \times Y} \cong H_X \otimes H_Y
$$
and the non-degeneracy of the intersection pairing (which even in the context of singular cohomology holds only after tensoring with $\BQ$, but in other situations such as Chow groups fails entirely). Therefore, for general cohomology theories, proving an identity between correspondences $\Gamma$ is a strictly stronger statement than the analogous identity between the operators $\Phi_\Gamma$. 

\medskip

\subsection{}

In the present Subsection, we will consider equivariant cohomology with respect to the torus $T = \BC^* \times \BC^*$. In particular, we have
$$
H^*_T(\text{point}) = \BZ[t_1,t_2]
$$
where $t_1,t_2$ are the standard coordinates on the Lie algebra of $T$. Recall the moduli space of framed sheaves of Subsection \ref{sub:framed}, and consider its equivariant cohomology
$$
H_{\CM_{\BA^2}} = H^*_T(\CM_{\BA^2}) = \bigoplus_{n = 0}^{\infty} H^*_T(\CM_{\BA^2,n})
$$
Recall the (natural analogue in the context of moduli spaces of framed sheaves of the) correspondence $\fZ$ of \eqref{eqn:simple}. Define operators
$$
e_n,f_n,h_n : H_{\CM_{\BA^2}} \rightarrow H_{\CM_{\BA^2}}
$$
for all $n \geq 0$ by the following formulas
\begin{equation}
\label{eqn:e yang}
e_n = t_1t_2 \cdot \pi_{+*}(\ell^n \cdot \pi_-^*)  \qquad \ 
\end{equation}
\begin{equation}
\label{eqn:f yang}
f_n = t_1t_2 \cdot \pi_{-*}(\ell^n (-1)^r \cdot \pi_+^*) 
\end{equation}
where  $\ell = c_1(\CL)$ is the first Chern class of the tautological line bundle \eqref{eqn:tautological}, and
\begin{equation}
\label{eqn:h yang}
h(z) = 1 + \sum_{n=0}^{\infty} \frac {h_n}{z^{n+1}} =  \text{cup product with } \frac {c(\CU_{\circ}, z+t)}{c(\CU_{\circ},z)}
\end{equation}
where $\CU_{\circ}$ denotes the restricted universal sheaf \eqref{eqn:restricted universal}, and $t=t_1+t_2$.

\medskip

\begin{remark}

Although the prefactor $t_1t_2$ in \eqref{eqn:e yang} and \eqref{eqn:f yang} might seem strange, it was chosen so that $e_n$ equals the composition
$$
H_{\CM_{\BA^2}} \xrightarrow{(\pi_+ \times \pi_{\BA^2})_*(\ell^n \cdot \pi_-^*)} H_{\CM_{\BA^2} \times \BA^2} \xrightarrow{|_{\circ}} H_{\CM_{\BA^2}}
$$
where $|_\circ$ denotes restriction to the origin of $\BA^2$ (and similarly  for the $f$ operators). This point of view will lead into the case of general surfaces $S$ in the next Subsection.

\end{remark}

\medskip

\noindent The connection between the abstract Yangian of Definition \ref{def:yang} and the operators above is the following result (see \cite[Theorem 3.8]{N ext} for notation closer to ours).

\medskip

\begin{theorem}
\label{thm:yang}

(the version of \cite{varagnolo} for the Lie algebra $\hgl_1$) For any $r$, the operators \eqref{eqn:e yang}, \eqref{eqn:f yang}, \eqref{eqn:h yang} satisfy relations \eqref{eqn:rel 1 yang}--\eqref{eqn:rel 5 yang}, thus yielding an action
$$
Y_{t_1,t_2}(\hgl_1) \curvearrowright H_{\CM_{\BA^2}}
$$

\end{theorem}

\medskip

\subsection{} 
\label{sub:operators cohomology}

We will now consider the case of general smooth projective surfaces $S$, and the moduli spaces of stable sheaves studied in Subsection \ref{sub:basic moduli}. Let us fix $(r,c_1) \in \BN \times H^2(S,\BZ)$ satisfying Assumptions A and S (from \eqref{eqn:assumption a} and \eqref{eqn:assumption s}, respectively) and consider the singular cohomology with integer coefficients
$$
H_{\CM \times S^k} = \bigoplus_{c_2 \in \BZ} H_{\CM_{(r,c_1,c_2)} \times S^k}
$$
for any $k \geq 0$. Due to Bogomolov's inequality \eqref{eqn:bogomolov}, the direct summands above are zero for $c_2$ small enough. Recall the variety $\fZ$ of \eqref{eqn:simple}, and define the operators
\begin{equation}
\label{eqn:e coh}
e_n = (\pi_+ \times \pi_S)_*(\ell^n \cdot \pi_-^*) : H_{\CM} \rightarrow H_{\CM \times S}  \qquad \ 
\end{equation}
\begin{equation}
\label{eqn:f coh}
f_n = (\pi_- \times \pi_S)_*(\ell^n (-1)^r\cdot \pi_+^*) : H_{\CM} \rightarrow H_{\CM \times S} 
\end{equation}
for all $n \geq 0$, where $\ell = c_1(\CL)$. We will also encounter the composition
\begin{equation}
\label{eqn:h coh series}
h(z) = 1 +  \sum_{n=0}^{\infty} \frac {h_n}{z^{n+1}} : H_{\CM} \xrightarrow{\text{pull-back}} H_{\CM \times S} \xrightarrow{\smile \frac {c(\CU,z+t)}{c(\CU,z)}} H_{\CM \times S}
\end{equation}
where we let $t=c_1(\CK_S)$. Combine the operators \eqref{eqn:e coh} and \eqref{eqn:f coh} into power series
\begin{equation}
\label{eqn:e coh series}
e(z) = \sum_{n=0}^{\infty} \frac {e_n}{z^{n+1}} = (\pi_+ \times \pi_S)_*\left(\frac 1{z-\ell} \cdot \pi_-^* \right)
\end{equation}
\begin{equation}
\label{eqn:f coh series}
f(z) = \sum_{n=0}^{\infty} \frac {f_n}{z^{n+1}} = (\pi_- \times \pi_S)_*\left(\frac {(-1)^r}{z-\ell} \cdot \pi_+^* \right)
\end{equation}

\medskip

\subsection{}
\label{sub:compositions}

We are now ready for the main result of the present Subsection, which describes the commutation relations between the operators $e_n, f_n, h_n$. We will find that these relations look best when expressed in terms of the series $e(z), f(z), h(z)$. This will involve the square $\textcolor{red}{S} \times \textcolor{blue}{S}$ of our surface; we choose to color the two factors of $S$ in red and blue, to make it easier to keep track of them in our formulas. Let
\begin{equation}
\label{eqn:def zeta coh}
\zeta^\coh(x) = 1 + \frac {[\Delta]}{x(x+t)} \in H_{\textcolor{red}{S} \times \textcolor{blue}{S}}(x)
\end{equation}
where $[\Delta] \in H_{S \times S}$ is the class of the diagonal, and 
\begin{equation}
\label{eqn:def zeta coh tilde} 
\tzeta^{\coh}_\pm(x) = \zeta(x)(x\pm \textcolor{red}{t})(x\mp \textcolor{blue}{t}) \in \frac {H_{\textcolor{red}{S} \times \textcolor{blue}{S}}[x]}x 
\end{equation}
where $\textcolor{red}{t}$ and $\textcolor{blue}{t}$ denote the pull-backs of the class $t=c_1(\CK_S)$ to $H_{\textcolor{red}{S} \times \textcolor{blue}{S}}$ via the first and second projection, respectively. In \eqref{eqn:def zeta coh}, the symbol $t$ in the denominator denotes either $\textcolor{red}{t}$ or $\textcolor{blue}{t}$; it does not matter which due to the presence of $[\Delta]$.

\medskip

\noindent Let us consider any operators $x,y : H_{\CM} \rightarrow H_{\CM \times S}$, for example \eqref{eqn:e coh}, \eqref{eqn:f coh} and the coefficients of \eqref{eqn:h coh series}. We may form the following two compositions 
\begin{align}
&\textcolor{red}{x} \textcolor{blue}{y} : H_{\CM} \xrightarrow{\textcolor{blue}{y}} H_{\CM \times \textcolor{blue}{S}} \xrightarrow{\textcolor{red}{x} \boxtimes \text{Id}_{\textcolor{blue}{S}}} H_{\CM\times \textcolor{red}{S} \times \textcolor{blue}{S}} \label{eqn:composition xy} \\
& \textcolor{blue}{y} \textcolor{red}{x} : H_{\CM} \xrightarrow{\textcolor{red}{x}} H_{\CM \times \textcolor{red}{S}} \xrightarrow{\textcolor{blue}{y} \boxtimes \text{Id}_{\textcolor{red}{S}}} H_{\CM\times \textcolor{red}{S} \times \textcolor{blue}{S}} \label{eqn:composition yx}
\end{align}
The colors of the operators will always match the color of the factor of $\textcolor{red}{S} \times \textcolor{blue}{S}$ in which they operate. We will write
$$
[\textcolor{red}{x}, \textcolor{blue}{y}] = \textcolor{red}{x} \textcolor{blue}{y} - \textcolor{blue}{y} \textcolor{red}{x}
$$
and we let $\Delta: S \hookrightarrow \textcolor{red}{S} \times \textcolor{blue}{S}$ denote the diagonal.

\medskip

\begin{theorem}
\label{thm:coh}

We have the following equalities of operators $H_{\CM} \rightarrow H_{\CM \times \textcolor{red}{S} \times \textcolor{blue}{S}}$
\begin{equation}
\label{eqn:rel 1 coh}
\left[ \textcolor{red}{e(z)} \textcolor{blue}{e(w)} \tzeta^{\ecoh}_{-}(w-z) \right]_{z^{<0},w^{<0}} = \left[ \textcolor{blue}{e(w)}\textcolor{red}{e(z)}  \tzeta^{\ecoh}_{+}(z-w) \right]_{z^{<0},w^{<0}}
\end{equation}
\begin{equation}
\label{eqn:rel 2 coh}
\left[ \textcolor{blue}{f(w)} \textcolor{red}{f(z)}  \tzeta^{\ecoh}_{-}(w-z) \right]_{z^{<0},w^{<0}} = \left[ \textcolor{red}{f(z)} \textcolor{blue}{f(w)} \tzeta^{\ecoh}_{+}(z-w) \right]_{z^{<0},w^{<0}}
\end{equation}
\begin{equation}
\label{eqn:rel 3 coh}
\textcolor{red}{h(z)} \textcolor{blue}{e(w)} =  \left[ \textcolor{blue}{e(w)} \textcolor{red}{h(z)} \frac {\zeta^{\ecoh}(z-w)}{\zeta^{\ecoh}(w-z)} \right]_{z \gg w|z^{\leq 0},w^{<0}}
\end{equation}
\begin{equation}
\label{eqn:rel 4 coh}
\textcolor{blue}{f(w)} \textcolor{red}{h(z)} = \left[ \textcolor{red}{h(z)} \textcolor{blue}{f(w)} \frac {\zeta^{\ecoh}(z-w)}{\zeta^{\ecoh}(w-z)} \right]_{z \gg w|z^{\leq 0},w^{<0}}
\end{equation}
\begin{equation}
\label{eqn:rel 5 coh}
\left [\textcolor{red}{f(z)}, \textcolor{blue}{e(w)} \right] = \frac 1{t} \cdot \Delta_* \left( \frac {h(z)-h(w)}{z-w} \right)
\end{equation}
as well as $[\textcolor{red}{h(z)}, \textcolor{blue}{h(w)}] = 0$. The meaning of all the square brackets is exactly the same as in Definition \ref{def:yang}. The class $t$ in the denominator of \eqref{eqn:rel 5 coh} does not pose an issue because all $h_n$ with $n > 0$ are multiples of $t$, as is clear from \eqref{eqn:h coh series}. 

\end{theorem}

\medskip

\noindent The proof of the Theorem above closely follows that of the upcoming Theorem \ref{thm:k-theory} in $K$-theory, and the former can actually be deduced from the latter via the Chern character isomorphism. We will therefore leave Theorem \ref{thm:coh} to the reader. 

\medskip

\subsection{}

In the remainder of the present Section, we will compare Theorem \ref{thm:yang} with Theorem \ref{thm:coh}. On one hand, the former is the version of the latter when $S = \BA^2$, stable sheaves are replaced by framed sheaves, and one works with equivariant cohomology (specifically, the operators \eqref{eqn:e yang}, \eqref{eqn:f yang}, \eqref{eqn:h yang} are obtained by composing the operators \eqref{eqn:e coh}, \eqref{eqn:f coh}, \eqref{eqn:h coh series} with the restriction map to the origin
$$
H_{\CM \times \BA^2} \xrightarrow{|_\circ} H_{\CM}
$$
which is an isomorphism). However, the more general Theorem \ref{thm:coh} can also be construed as saying that ``the Yangian acts on $H_{\CM}$" for a general surface $S$, although care must be taken to properly formulate what this means. This will be the goal of the remainder of the present Section, and we will start by showing how to unravel relations \eqref{eqn:rel 1 coh}--\eqref{eqn:rel 5 coh}. We start with \eqref{eqn:rel 5 coh}, whose $z^{-n-1}w^{-m-1}$ coefficient reads
\begin{equation}
\label{eqn:rel 5 coh explicit}
\left[ \textcolor{red}{f_{n}}, \textcolor{blue}{e_{m}} \right] = \Delta_*\left(- \frac {h_{m+n}}t \right)
\end{equation}
as an equality of operators $H_{\CM} \rightarrow H_{\CM \times \textcolor{red}{S} \times \textcolor{blue}{S}}$ (compare with \eqref{eqn:rel 5 yang explicit}). 

\medskip

\begin{remark}

If one works with cohomology with rational coefficients instead of integer coefficients (i.e. assume that $H_X$ is defined as $H^*(X,\BQ)$ in this Remark only), one may use the K\"unneth decomposition
$$
H_{\CM \times S} \cong H_{\CM} \otimes H_S
$$
with respect to which the operators $e_m$, $f_n$, $h_{m+n}$ decompose as
$$
e_m = \sum_{i} E_{(\gamma_i)} \otimes \gamma^i, \quad f_n = \sum_{i} F_{(\gamma_i)} \otimes \gamma^i, \quad \frac {h_{n+m}}t = \sum_{i} H_{(\gamma_i)} \otimes \gamma^i
$$
where $\{\gamma_i,\gamma^i\}$ run over fixed dual bases of $H_S$, and $E_{(\gamma_i)}, F_{(\gamma_i)}, H_{(\gamma_i)} : H_{\CM} \rightarrow H_{\CM}$. In this case, \eqref{eqn:rel 5 coh explicit} is equivalent to the following equality of operators $H_{\CM} \rightarrow H_{\CM}$
\begin{equation}
\label{eqn:rel 5 coh kunneth}
\left [F_{(\gamma_i)}, E_{(\gamma_j)} \right] = - H_{(\gamma_i\gamma_j)}, \quad \forall  i,j
\end{equation}
(in the right-hand side, we define $H_{(\gamma_i\gamma_j)}$ as the same linear combination of $H_{(\gamma_k)}$'s as $\gamma_i\gamma_j$ is a linear combination of $\gamma_k$'s). While formulas such as \eqref{eqn:rel 5 coh kunneth} are used more commonly in the literature, we prefer formula \eqref{eqn:rel 5 coh explicit} because its holds for cohomology theories without K\"unneth decompositions, such as Chow groups.

\end{remark}

\medskip

\subsection{}

Let us now unravel relations \eqref{eqn:rel 3 coh} and \eqref{eqn:rel 4 coh}. We start by noting that 
$$
\frac {\zeta^{\coh}(z-w)}{\zeta^{\coh}(w-z)} = \frac {\displaystyle 1+\frac {[\Delta]}{(z-w)(z-w+t)}}{\displaystyle 1+\frac {[\Delta]}{(z-w)(z-w-t)}} = 1 + \sum_{a =3}^{\infty} \sum_{b=0}^{a-2} \Delta_*(\gamma_{ab}) \frac {w^b}{z^a}
$$
where $\gamma_{ab} \in H_S$ is the same polynomial expression in the classes
\begin{align*}
&t_1+t_2 = c_1(\Omega^1_S) \\
&t_1t_2 = c_2(\Omega^1_S)
\end{align*}
as the one we already encountered in \eqref{eqn:expansion yangian}. The Chern classes above arise because $t = t_1+t_2$ and because of the identity $[\Delta]  [\Delta] = t_1t_2[\Delta]$ in $H_{S \times S}$. Thus, relations \eqref{eqn:rel 3 coh} and \eqref{eqn:rel 4 coh} take the following form, for all $m,n \geq 0$
\begin{equation}
\label{eqn:rel 3 coh explicit}
[\textcolor{red}{h_{n}}, \textcolor{blue}{e_{m}}] = \Delta_* \left(  \sum_{a =3}^{\infty} \sum_{b=0}^{a-2} \gamma_{ab} \underbrace{\textcolor{blue}{e_{m+b}} \textcolor{red}{h_{n-a}} \Big|_\Delta}_{\text{an operator } H_{\CM} \rightarrow H_{\CM \times \textcolor{red}{S} \times \textcolor{blue}{S}} \xrightarrow{|_\Delta} H_{\CM \times S}} \right)
\end{equation}
\begin{equation}
\label{eqn:rel 4 coh explicit}
[\textcolor{blue}{f_{m}}, \textcolor{red}{h_{n}}] = \Delta_* \left(  \sum_{a =3}^{\infty} \sum_{b=0}^{a-2} \gamma_{ab} \underbrace{\textcolor{red}{h_{n-a}}\textcolor{blue}{f_{m+b}}  \Big|_\Delta}_{\text{an operator } H_{\CM} \rightarrow H_{\CM \times \textcolor{red}{S} \times \textcolor{blue}{S}} \xrightarrow{|_\Delta} H_{\CM \times S}} \right)
\end{equation}
(in the right-hand sides above, we make the convention that $h_{-1}$ is the pull-back map and $h_{-2} = h_{-3} = \dots = 0$). Compare the formulas above with \eqref{eqn:rel 3 yang explicit}--\eqref{eqn:rel 4 yang explicit}. 

\medskip

\noindent As for formulas \eqref{eqn:rel 1 coh}--\eqref{eqn:rel 2 coh}, let us first observe that
\begin{align*}
&\tzeta^{\coh}_{+}(z-w) = (z-w)^2 + (z-w)(\textcolor{red}{t} - \textcolor{blue}{t}) - \textcolor{red}{t} \textcolor{blue}{t} + [\Delta] - \frac{[\Delta]t}{z-w} \\
&\tzeta^{\coh}_{-}(w-z) = (z-w)^2 + (z-w)(\textcolor{red}{t} - \textcolor{blue}{t}) - \textcolor{red}{t} \textcolor{blue}{t} + [\Delta] + \frac{[\Delta]t}{z-w}
\end{align*}
Therefore, \eqref{eqn:rel 1 coh} and \eqref{eqn:rel 2 coh} yield the following formulas
\begin{equation}
\label{eqn:rel 1 coh explicit}
[\textcolor{red}{e_{n+3}},\textcolor{blue}{e_m}] - 3 [\textcolor{red}{e_{n+2}},\textcolor{blue}{e_{m+1}}] + 3 [\textcolor{red}{e_{n+1}},\textcolor{blue}{e_{m+2}}] - [\textcolor{red}{e_{n}},\textcolor{blue}{e_{m+3}}] +  
\end{equation}
$$
 + (\textcolor{red}{t} - \textcolor{blue}{t}) ([\textcolor{red}{e_{n+2}},\textcolor{blue}{e_{m}}] - 2[\textcolor{red}{e_{n+1}},\textcolor{blue}{e_{m+1}}] + [\textcolor{red}{e_{n}},\textcolor{blue}{e_{m+2}}]) - \textcolor{red}{t} \textcolor{blue}{t} ([\textcolor{red}{e_{n+1}},\textcolor{blue}e_m]-[\textcolor{red}{e_n},\textcolor{blue}{e_{m+1}}]) + 
$$
$$
+ \Delta_*\left(e_{n+1} e_m \Big|_\Delta - e_m e_{n+1} \Big|_\Delta - e_n e_{m+1} \Big|_\Delta + e_{m+1}e_n \Big|_\Delta +  te_ne_m \Big|_\Delta+ te_me_n\Big|_\Delta \right) = 0
$$
and 
\begin{equation}
\label{eqn:rel 2 coh explicit}
[\textcolor{blue}{f_m},\textcolor{red}{f_{n+3}}] - 3 [\textcolor{blue}{f_{m+1}},\textcolor{red}{f_{n+2}}] + 3 [\textcolor{blue}{f_{m+2}},\textcolor{red}{f_{n+1}}] - [\textcolor{blue}{f_{m+3}},\textcolor{red}{f_{n}}] +  
\end{equation}
$$
 + (\textcolor{red}{t} - \textcolor{blue}{t}) ([\textcolor{blue}{f_{m}},\textcolor{red}{f_{n+2}}] - 2[\textcolor{blue}{f_{m+1}},\textcolor{red}{f_{n+1}}] + [\textcolor{blue}{f_{m+2}},\textcolor{red}{f_{n}}]) - \textcolor{red}{t}  \textcolor{blue}{t} ([\textcolor{blue}{f_m},\textcolor{red}{f_{n+1}}]-[\textcolor{blue}{f_{m+1}},\textcolor{red}{f_n}]) + 
$$
$$
+ \Delta_*\left(f_m f_{n+1}  \Big|_\Delta - f_{n+1} f_m  \Big|_\Delta - f_{m+1}f_n  \Big|_\Delta + f_nf_{m+1} \Big|_\Delta +  tf_mf_n \Big|_\Delta+ tf_nf_m\Big|_\Delta \right) = 0
$$
for all $m,n \geq 0$.

\medskip

\subsection{}
\label{sub:diagonal 1}

We have phrased relations \eqref{eqn:rel 5 coh explicit}, \eqref{eqn:rel 3 coh explicit}, \eqref{eqn:rel 4 coh explicit}, \eqref{eqn:rel 1 coh explicit}, \eqref{eqn:rel 2 coh explicit} as equalities of compositions of operators, but it is better if one regards them as equalities of correspondences in $\CM \times \CM \times S \times S$, in the language of Subsection \ref{sub:correspondences}. For example, 
$$
\textcolor{red}{f_n}\textcolor{blue}{e_m} \quad \text{ is given by the class } \quad \ell_1^n \ell_2^m \quad \text{ on } \quad \Big\{\CF \supset_{\textcolor{red}{x}} \widetilde{\CF} \subset_{\textcolor{blue}{y}} \CF'\Big\}
$$
$$
\textcolor{blue}{e_m}\textcolor{red}{f_n} \quad \text{is given by the class} \quad {\ell_1'}^m {\ell_2'}^n \quad \text{on} \quad \Big\{\CF \subset_{\textcolor{blue}{y}} \widetilde{\CF}' \supset_{\textcolor{red}{x}} \CF'\Big\}
$$
where $\ell_1 = c_1(\CF_x/\widetilde{\CF}_x)$, $\ell_2 = c_1(\CF'_y/\widetilde{\CF}_y)$, $\ell_1'=c_1(\widetilde{\CF}'_y/\CF_y)$, $\ell_2'=c_1(\widetilde{\CF}'_x/\CF'_x)$. One regards the cohomology classes above as correspondences by pushing them forward under the forgetful maps
$$
\Big\{\CF \supset_{\textcolor{red}{x}} \widetilde{\CF} \subset_{\textcolor{blue}{y}} \CF'\Big\}, \Big\{\CF \subset_{\textcolor{blue}{y}} \widetilde{\CF}' \supset_{\textcolor{red}{x}} \CF'\Big\} \rightarrow \CM \times \CM \times \textcolor{red}{S} \times \textcolor{blue}{S}
$$
that only remember $(\CF,\CF',\textcolor{red}{x},\textcolor{blue}{y})$. These correspondences are actually equal on the open locus $\textcolor{red}{x} \neq \textcolor{blue}{y}$, because of the mutually inverse isomorphisms
\begin{align*}
&\Big\{\CF \supset_{\textcolor{red}{x}} \widetilde{\CF} \subset_{\textcolor{blue}{y}} \CF'\Big\} \longrightarrow \Big\{\CF \subset_{\textcolor{blue}{y}} \widetilde{\CF}' \supset_{\textcolor{red}{x}} \CF'\Big\}, \qquad \widetilde{\CF}' = \CF \oplus_{\widetilde{\CF}} \CF' \\
&\Big\{\CF \subset_{\textcolor{blue}{y}} \widetilde{\CF}' \supset_{\textcolor{red}{x}} \CF'\Big\} \longrightarrow \Big\{\CF \supset_{\textcolor{red}{x}} \widetilde{\CF} \subset_{\textcolor{blue}{y}} \CF'\Big\}, \qquad \widetilde{\CF} = \CF \cap \CF' \text{ inside } \widetilde{\CF}' 
\end{align*}
which identify $\ell_1 = \ell_2'$ and $\ell_2 = \ell_1'$ on the locus $\{\textcolor{red}{x} \neq \textcolor{blue}{y}\}$. Therefore, we infer that
$$
[\textcolor{red}{f_n}, \textcolor{blue}{e_m}] \Big|_{\{\textcolor{red}{x} \neq \textcolor{blue}{y}\}} = 0
$$
which by the excision long exact sequence in cohomology yields
$$
[\textcolor{red}{f_n}, \textcolor{blue}{e_m}] = \Delta_*(c)
$$
for some operator $c : H_{\CM} \rightarrow H_{\CM \times S}$. The content of \eqref{eqn:rel 5 coh explicit} is that $c$ is equal to the operator of cup product with a certain specific cohomology class. 

\medskip

\subsection{}
\label{sub:diagonal 2}

The geometric argument presented in the previous Subsection actually holds for any operators $a,b \in \{e_n,f_n,h_n\}_{n \geq 0}$. Thus, the commutator $[\textcolor{red}{a}, \textcolor{blue}{b}]$ vanishes away from the diagonal of $\textcolor{red}{S} \times \textcolor{blue}{S}$, so the excision long exact sequence implies that
$$
[\textcolor{red}{a}, \textcolor{blue}{b}] = \Delta_*(c)
$$
for some $c : H_{\CM} \rightarrow H_{\CM \times S}$. Because $\Delta_*$ is injective (projection on one of the factors is a left inverse of $\Delta$), such a $c$ is unique, and we will therefore denote it by
\begin{equation}
\label{eqn:reduced}
[\textcolor{red}{a}, \textcolor{blue}{b}]_{\text{red}} =: c
\end{equation}
In this case, we will say that $a$ and $b$ \textbf{have diagonal commutator}. We are now ready to define the Yangian action on $H_{\CM}$. Consider the ring homomorphism
$$
\phi : \BZ[t_1,t_2]^{\sym} \rightarrow H_S, \qquad \phi(t_1+t_2) = c_1(\Omega_S^1), \qquad \phi(t_1t_2) = c_2(\Omega_S^1)
$$
and let $\CM \times S \xrightarrow{\pi} \CM$ and $\CM \times S \xrightarrow{\rho} S$ denote the standard projections.

\medskip

\begin{definition}
\label{def:action}

An action $Y_{t_1,t_2}(\hgl_1) \curvearrowright H_{\CM}$ is an abelian group homomorphism
$$
Y_{t_1,t_2}(\hgl_1) \xrightarrow{\Phi} \emph{Hom}(H_{\CM}, H_{\CM \times S})
$$
satisfying the following properties for all $x,y \in Y_{t_1,t_2}(\hgl_1)$ and $\gamma \in \BZ[t_1,t_2]^{\emph{sym}}$

\medskip

\begin{itemize}[leftmargin=*]

\item \textbf{unit:} 
\begin{equation}
\label{eqn:unit}
\Phi(1) = \Big( H_{\CM} \xrightarrow{\pi^*} H_{\CM \times S} \Big)
\end{equation}

\medskip

\item $\BZ[t_1,t_2]^{\esym}$-\textbf{linearity:} 
\begin{equation}
\label{eqn:linearity}
\Phi(\gamma x) = \Big( H_{\CM} \xrightarrow{\Phi(x)} H_{\CM \times S} \xrightarrow{\smile \rho^*(\phi(\gamma))} H_{\CM \times S} \Big) 
\end{equation}

\medskip

\item \textbf{multiplicativity:}  
\begin{equation}
\label{eqn:multiplicativity}
\Phi(xy) = \Big( H_{\CM} \xrightarrow{\Phi(y)} H_{\CM \times \textcolor{blue}{S}} \xrightarrow{\Phi(x) \boxtimes \emph{Id}_{\textcolor{blue}{S}}} H_{\CM\times \textcolor{red}{S} \times \textcolor{blue}{S}} \xrightarrow{|_\Delta} H_{\CM \times S}\Big)
\end{equation}

\medskip

\item \textbf{commutator:} $\Phi(x)$ and $\Phi(y)$ have diagonal commutator, and
\begin{equation}
\label{eqn:commutator}
[\Phi(x),\Phi(y)]_{\emph{red}} = \Phi \left(\frac {[x,y]}{t_1t_2}\right)
\end{equation}

\end{itemize}

\noindent The right-hand side of \eqref{eqn:commutator} is well-defined because of \eqref{eqn:divisible}.

\end{definition}

\medskip

\subsection{} Comparing relations \eqref{eqn:rel 1 yang}--\eqref{eqn:rel 5 yang} with \eqref{eqn:rel 1 coh}--\eqref{eqn:rel 5 coh}, one can translate Theorem \ref{thm:coh} into the existence of an action
\begin{equation}
\label{eqn:action yangian proof}
Y_{t_1,t_2}(\hgl_1) \curvearrowright H_{\CM}
\end{equation}
as stated in Theorem \ref{thm:coh intro}. The main novelty behind the geometric notion of action in Definition \ref{def:action} is that the commutator and product of operators are objects of a different nature: the former is an operator $H_{\CM} \rightarrow H_{\CM \times \textcolor{red}{S} \times \textcolor{blue}{S}}$ and the latter is an operator $H_{\CM} \rightarrow H_{\CM \times S}$. However, if we restrict to operators $a,b : H_{\CM} \rightarrow H_{\CM \times S}$ with diagonal commutator, then
$$
[a,b]_{\text{red}} \quad \text{and} \quad ab|_\Delta \quad \text{are both operators} \quad H_{\CM} \rightarrow H_{\CM \times S}
$$
Moreover, these operators satisfy the following natural analogues of associativity
\begin{equation}
\label{eqn:associativity}
\left(ab|_\Delta c\right)\Big|_\Delta= a \left(bc|_\Delta\right)\Big|_\Delta
\end{equation}
the Leibniz rule 
\begin{equation}
\label{eqn:associativity}
\left[a,bc|_\Delta \right]_{\text{red}} = [a,b]_{\text{red}}c\Big|_\Delta + b [a,c]_{\text{red}}\Big|_\Delta
\end{equation}
and the Jacobi identity
\begin{equation}
\label{eqn:associativity}
\left[[a,b]_{\text{red}},c\right]_{\text{red}} + \left[[b,c]_{\text{red}},a\right]_{\text{red}} + \left[[c,a]_{\text{red}},b\right]_{\text{red}} = 0
\end{equation}
for any $a,b,c : H_{\CM} \rightarrow H_{\CM \times S}$ for which all $[\cdot,\cdot]_{\text{red}}$ above are defined.

\bigskip

\section{$K$-theory}
\label{sec:k-theory}

\medskip

\subsection{} 
\label{sub:tor} 

In the present Section, we will present trigonometric/$K$-theoretic versions of the rational/cohomological constructions in the previous Section. We will work over the ring of Laurent polynomials in two variables $q_1$ and $q_2$, which should be interpreted as the exponentials of the parameters $t_1$ and $t_2$ from Subsection \ref{sub:yang}. Let $q = q_1q_2$ and consider the rational functions 
\begin{equation}
\label{eqn:def zeta tor}
\zeta^{\BC^*}(x) = \frac {(1-xq_1)(1-xq_2)}{(1-x)(1-xq)}
\end{equation}
\begin{equation}
\label{eqn:def zeta tor tilde}
\tzeta^{\BC^*}(x) = \zeta^{\BC^*}(x)(1-xq)(1-x^{-1}q) = \frac {(1-xq_1)(1-xq_2)(1-x^{-1}q)}{1-x}
\end{equation}
by analogy with \eqref{eqn:def zeta yang} and \eqref{eqn:def zeta yang tilde}. 

\medskip

\begin{definition}
\label{def:tor}

Quantum toroidal $\fgl_1$ is the algebra
$$
U_{q_1,q_2}(\ddot{\fgl}_1) = \BZ[q^{\pm 1}_1,q^{\pm 1}_2] \Big\langle e_n, f_n, h^\pm_m \Big \rangle_{m \geq 0, n \in \BZ} \Big / \text{relations \eqref{eqn:rel 1 tor}--\eqref{eqn:rel 5 tor}}
$$
The defining relations are best written in terms of the generating series
$$
e(z) = \sum_{n=-\infty}^{\infty} \frac {e_n}{z^n}, \qquad f(z) = \sum_{n=-\infty}^{\infty} \frac {f_n}{z^n}, \qquad h^\pm(z) = \sum_{m=0}^{\infty} \frac {h^\pm_m}{z^{\pm m}}
$$
and take the form
\begin{equation}
\label{eqn:rel 1 tor}
e(z)e(w) \tzeta^{\BC^*}\left(\frac wz \right) = e(w)e(z) \tzeta^{\BC^*}\left(\frac zw \right)
\end{equation}
\begin{equation}
\label{eqn:rel 2 tor}
f(w) f(z) \tzeta^{\BC^*}\left(\frac wz \right) = f(z)f(w) \tzeta^{\BC^*}\left(\frac zw \right)
\end{equation}
\begin{equation}
\label{eqn:rel 3 tor}
h^\pm(z)e(w) =  e(w)h^\pm(z) \frac {\zeta^{\BC^*}\left(\frac zw \right)}{\zeta^{\BC^*}\left(\frac wz \right)}
\end{equation}
\begin{equation}
\label{eqn:rel 4 tor}
f(w) h^\pm(z) = h^\pm(z) f(w) \frac {\zeta^{\BC^*}\left(\frac zw \right)}{\zeta^{\BC^*}\left(\frac wz \right)}
\end{equation}
\begin{equation}
\label{eqn:rel 5 tor}
\left [f(z), e(w) \right] = \frac {(1-q_1)(1-q_2)}{1-q} \cdot \delta \left(\frac zw\right) (h^+(z)-h^-(w))
\end{equation}
(where $\delta(x) = \sum_{n=-\infty}^{\infty} x^n$), as well as $[h^\pm(z),h^{\pm'}(w)] = 0$ for all $\pm, \pm' \in \{+,-\}$.

\end{definition}

\medskip

\noindent In formulas \eqref{eqn:rel 1 tor}--\eqref{eqn:rel 2 tor}, we cancel out the factor $z-w$ from the denominator and equate the coefficients of all $\{z^aw^b\}_{a,b\in \BZ}$ in the left and right-hand sides. In the next two formulas \eqref{eqn:rel 3 tor}--\eqref{eqn:rel 4 tor}, we expand the rational function in the RHS in non-negative powers of $w^{\pm 1}/z^{\pm 1}$, and equate the coefficients of all $\{z^aw^b\}_{\pm a \in \BZ_{\leq 0}, b \in \BZ}$ in the left and right-hand sides. Finally, in the last formula \eqref{eqn:rel 5 tor}, we simply equate the coefficients of all $\{z^aw^b\}_{a,b \in \BZ}$ in the left and right-hand sides.

\medskip

\begin{remark}
\label{rem:extra relation tor}

Definition \ref{def:tor} is the original one of Ding-Iohara (\cite{di}). Similarly with Remark \ref{rem:extra relation yang}, it is customary to impose the additional cubic relations
\begin{equation}
\label{eqn:cubic tor}
\sum_{\sigma \in S_3} [e_{n_{\sigma(1)}},[e_{n_{\sigma(2)}-1},e_{n_{\sigma(3)}+1}]]  =  \sum_{\sigma \in S_3} [f_{n_{\sigma(1)}},[f_{n_{\sigma(2)}-1},f_{n_{\sigma(3)}+1}]]  =  0
\end{equation}
in $\UU$, for all $n_1,n_2,n_3 \in \BZ$. This was the approach of Miki (\cite{miki}), who imposed the particular case of the cubic relations \eqref{eqn:cubic tor} for $n_1=n_2=n_3$ (the more general versions of relations \eqref{eqn:cubic tor} appeared in \cite{tsymbaliuk}, and they actually follow from the particular cases studied by Miki and \eqref{eqn:rel 1 tor}--\eqref{eqn:rel 5 tor}). While we will not consider the cubic relations in the present paper for brevity, we observe that they do indeed hold in all the modules of geometric nature considered in the present paper, namely the $K$-theory groups of moduli spaces of stable/framed sheaves (as shown in \cite{N hecke}).

\end{remark}

\medskip

\subsection{}

Let us unpack the relations in Definition \ref{def:tor}. The easiest is \eqref{eqn:rel 5 tor}, which states
\begin{equation}
\label{eqn:rel 5 tor explicit}
[f_n,e_m] = \frac {(1-q_1)(1-q_2)}{1-q} \cdot \begin{cases} h^+_{m+n} &\text{if }m+n > 0 \\ h^+_0 - h^-_0 &\text{if }m+n =0 \\ -h^-_{-m-n}&\text{if }m+n < 0 \end{cases}
\end{equation}
for all $m,n \in \BZ$. To express relations \eqref{eqn:rel 3 tor} and \eqref{eqn:rel 4 tor}, we need to consider the power series expansion
\begin{equation}
\label{eqn:expansion toroidal}
\frac {\zeta^{\BC^*}\left(\frac zw\right)}{\zeta^{\BC^*}\left(\frac wz\right)} = \frac {(zq_1-w)(zq_2-w)(z-wq)}{(z-wq_1)(z-wq_2)(zq-w)} = 1 + \sum_{a=1}^{\infty} (1-q_1)(1-q_2) \gamma^\pm_{a} \frac {w^{\pm a}}{z^{\pm a}}
\end{equation}
for various Laurent polynomials $\gamma^\pm_{a}$ in $q=q_1q_2$ and $q_1+q_2$. The fact that we can extract a factor of $(1-q_1)(1-q_2)$ from the coefficients in the right-hand side of \eqref{eqn:expansion toroidal} will be very important for our geometric constructions in the later Subsections, and it is due to the fact that the left-hand side of \eqref{eqn:expansion toroidal} is equal to 1 when either $q_1=1$ or when $q_2=1$. Therefore, relations \eqref{eqn:rel 3 tor}--\eqref{eqn:rel 4 tor} read
\begin{align}
&[h^\pm_n,e_m] = (1-q_1)(1-q_2) \sum_{a=1}^{n} \gamma^\pm _{a} \cdot e_{m\pm a}h^\pm_{n-a} \label{eqn:rel 3 tor explicit} \\
&[f_m,h^\pm_n] = (1-q_1)(1-q_2) \sum_{a=1}^{n} \gamma^\pm _{a} \cdot h^\pm_{n-a}f_{m\pm a} \label{eqn:rel 4 tor explicit}
\end{align}
for all $m \in \BZ$, $n \geq 0$. Finally, to unpack relations \eqref{eqn:rel 1 tor}--\eqref{eqn:rel 2 tor}, we observe that
\begin{multline*}
\widetilde{\zeta}^{\BC^*}\left( \frac zw \right) = \frac {(zq_1-w)(zq_2-w)(z-wq)}{zw(z-w)} = \\ = \left(1-\frac {zq}w\right)\left(1-\frac {wq}z\right) - (1-q_1)(1-q_2)\frac {z-wq}{z-w}
\end{multline*}
Therefore, relations \eqref{eqn:rel 1 tor}--\eqref{eqn:rel 2 tor} respectively imply that for all $m,n \in \BZ$ we have
$$
[e_{n+3},e_m] - \left(q+1+q^{-1} \right)[e_{n+2},e_{m+1}]  + \left(q+1+q^{-1} \right)[e_{n+1},e_{m+2}] - [e_n,e_{m+3}] + 
$$
\begin{equation}
\label{eqn:rel 1 tor explicit}
+ (1-q_1)(1-q_2) \left(e_{n+2}e_{m+1} - \frac {e_{n+1}e_{m+2}}q - \frac {e_{m+1}e_{n+2}}q + e_{m+2}e_{n+1}\right) = 0
\end{equation}
$$
[f_m,f_{n+3}] - \left(q+1+q^{-1} \right)[f_{m+1},f_{n+2}]  + \left(q+1+q^{-1} \right)[f_{m+2},f_{n+1}] - [f_{m+3},f_n] + 
$$
\begin{equation}
\label{eqn:rel 2 tor explicit}
+ (1-q_1)(1-q_2) \left(f_{m+1}f_{n+2} - \frac {f_{m+2}f_{n+1}}q - \frac {f_{n+2}f_{m+1}}q + f_{n+1}f_{m+2}\right) = 0
\end{equation}
Although we will not include this explicitly in our notation in order to not overburden it, we note that future connections with geometry require us to make the following two modifications to quantum toroidal $\fgl_1$, by analogy with Remark \ref{rem:quotients}

\medskip

\begin{itemize}[leftmargin=*]

\item $\UU$ should be defined over $\BZ[q_1^{\pm 1},q_2^{\pm 1}]^{\text{sym}}$ instead of over $\BZ[q_1^{\pm 1},q_2^{\pm 1}]$

\medskip

\item for any $x,y \in \UU$,
\begin{equation}
\label{eqn:divisible toroidal}
[x,y] \ \text{should be divisible by} \ (1-q_1)(1-q_2)
\end{equation} 
The way to achieve this is to formally adjoin the symbols 
$$
\frac {[e_{n_1},\dots,[e_{n_{k-1}},[e_{n_k},e_{n_{k+1}}]] \dots]}{(1-q_1)^k(1-q_2)^k} \quad \text{and} \quad \frac {[f_{n_1},\dots,[f_{n_{k-1}},[f_{n_k},f_{n_{k+1}}]] \dots]}{(1-q_1)^k(1-q_2)^k}
$$
to $\UU$, for all $n_1,\dots,n_{k+1} \in \BZ$, and impose the natural Leibniz rules and Jacobi identities between these new symbols and the old generators $e_n,f_n,h^\pm_m$. 

\end{itemize}

\medskip

\subsection{}

Given a smooth algebraic variety $X$ over $\BC$, its (0-th) algebraic $K$-theory group is defined as
\begin{equation}
\label{eqn:k-theory ring}
K_X = \bigoplus_{\CV \text{ locally free sheaf on }X} \BZ [\CV] \Big/\text{ relation \eqref{eqn:rel k-theory}}
\end{equation}
where for any short exact sequence 
\begin{equation}
\label{eqn:ses k-theory}
0 \rightarrow \CA \rightarrow \CV \rightarrow \CB \rightarrow 0
\end{equation}
we impose the relation
\begin{equation}
\label{eqn:rel k-theory}
[\CV] = [\CA]+[\CB]
\end{equation}
Algebraic $K$-theory is a ring with respect to tensor product, and it possesses pull-backs (for general maps) and push-forwards (for proper maps of smooth varieties).

\medskip

\begin{remark}

One can regard $X$ as a manifold and replace algebraic $K$-theory by topological $K$-theory, simply by changing ``locally free sheaf" with ``vector bundle" in \eqref{eqn:k-theory ring}. All the results in the present Section hold in the topological setup as well.

\end{remark}

\medskip 

\noindent In $K$-theory, the role of the Chern classes \eqref{eqn:chern} is played by the exterior powers
$$
[\wedge^i\CV] \in K_X
$$
for any locally free sheaf $\CV$ on $X$. We will combine these  in the Laurent polynomial
$$
\wedge^\bullet\left(\frac {\CV}z\right) = \sum_{i=0}^{\text{rank }\CV} (-z)^{-i} [\wedge^i\CV] \in K_X[z^{-1}]
$$
For any short exact sequence \eqref{eqn:ses} with $\CV,\CW$ locally free of ranks $v,w$, we may set
$$
\wedge^\bullet\left(\frac {\CE}z\right) = \frac {\wedge^\bullet\left(\frac {\CV}z\right)}{\wedge^\bullet\left(\frac {\CW}z\right)}
$$
and use this to unambigiously define the exterior powers of $\CE$
$$
\wedge^\bullet\left(\frac {\CE}z\right) = \sum_{i=0}^{\infty} (-z)^{-i} [\wedge^i\CE] \in K_X[[z^{-1}]]
$$
We will also use the following notation for the exterior powers of $\CE^\vee$
$$
\wedge^\bullet\left(\frac z{\CE}\right) = \wedge^\bullet\left(z \CE^\vee\right)
$$
If a smooth variety $X$ is endowed with an action of a torus $T$, then we may define its equivariant $K$-theory
\begin{equation}
\label{eqn:equivariant k-theory ring}
K_T(X)
\end{equation}
by considering $T$-equivariant locally free sheaves in \eqref{eqn:k-theory ring} and $T$-equivariant short exact sequences in \eqref{eqn:ses k-theory}. Very importantly, $K_T(X)$ is a module over the ring
$$
K_T(\text{point}) = \text{Rep}_T
$$
Indeed, any one-dimensional representation $\chi$ of $T$ can be tensored with an arbitrary $T$-equivariant locally free sheaf $\CV$; the resulting $\chi \otimes \CV$ is isomorphic to $\CV$ as a locally free sheaf, but not as a $T$-equivariant locally free sheaf. The aforementioned module structure is therefore given by $\chi \cdot [\CV] = [\chi \otimes \CV]$. Finally, note that the formalism of correspondences from Subsection \ref{sub:correspondences} holds equally well in $K$-theory (be it equivariant or not) simply by replacing $H$ with $K$ everywhere.

\medskip

\subsection{} In the present Subsection, we will consider $K$-theory equivariant with respect to the torus $T = \BC^* \times \BC^*$. Therefore, we have
$$
K_T(\text{point}) = \BZ[q^{\pm 1}_1,q^{\pm 1}_2]
$$
where $q_1,q_2$ are the standard characters of the two factors of $T$. We recall the moduli space of framed sheaves of Subsection \ref{sub:framed}, and consider its equivariant $K$-theory
$$
K_{\CM_{\BA^2}} = K_T(\CM_{\BA^2}) = \bigoplus_{n = 0}^{\infty} K_T(\CM_{\BA^2,n})
$$
Recall the (natural analogue in the context of moduli spaces of framed sheaves of the) correspondence $\fZ$ of \eqref{eqn:simple}, and we will define operators
$$
e_n,f_n,h^\pm_m : K_{\CM_{\BA^2}} \rightarrow K_{\CM_{\BA^2}}
$$
For all $n \in \BZ$, we set 
\begin{align}
&e_n = (1-q_1)(1-q_2) \cdot \pi_{+*} \Big(L^n \cdot \pi_-^*\Big)  \label{eqn:e tor}  \\
&f_n = (1-q_1)(1-q_2) \cdot \pi_{-*} \Big(L^{n-r} (-1)^r \det \CU_\circ \cdot \pi_+^*\Big)  \label{eqn:f tor}
\end{align}
where $L = [\CL]$ denotes the $K$-theory class of the tautological line bundle \eqref{eqn:tautological}, and $\CU_{\circ}$ is the restricted universal sheaf \eqref{eqn:restricted universal}. Moreover, let $q=q_1q_2$ and consider
\begin{equation}
\label{eqn:h tor}
h^\pm(z) = \sum_{m=0}^{\infty} \frac {h_m^\pm}{z^{\pm m}} =  \text{multiplication by }  \wedge^\bullet\left(\frac {z(q-1)}{\CU_\circ}\right)
\end{equation}
where the right-hand side is expanded in powers of $z^{\mp 1}$. The connection between the quantum toroidal $\fgl_1$ of Definition \ref{def:tor} and the operators above is the following.

\medskip

\begin{theorem}
\label{thm:tor}

(\cite{ft, sv}) For any $r \in \BN$, the operators \eqref{eqn:e tor}--\eqref{eqn:h tor} satisfy relations \eqref{eqn:rel 1 tor}--\eqref{eqn:rel 5 tor}, thus yielding an action
$$
\UU \curvearrowright K_{\CM_{\BA^2}}
$$

\end{theorem}

\medskip

\subsection{} 
\label{sub:operators k-theory}

We will now consider the case of general projective surfaces $S$, and the moduli spaces of stable sheaves studied in Subsection \ref{sub:basic moduli}. Thus, we will fix $(r,c_1) \in \BN \times H^2(S,\BZ)$ satisfying Assumptions A and S (from \eqref{eqn:assumption a} and \eqref{eqn:assumption s}) and consider 
$$
K_{\CM \times S^k} =  \bigoplus_{c_2 \in \BZ} K_{\CM_{(r,c_1,c_2)} \times S^k}
$$
for any $k \geq 0$. Due to Bogomolov's inequality \eqref{eqn:bogomolov}, the direct summands above are zero for $c_2$ small enough. Recall the variety $\fZ$ of \eqref{eqn:simple}, and define the operators
\begin{align}
&e_n = (\pi_+ \times \pi_S)_* \Big(L^n \cdot \pi_-^* \Big) \qquad \qquad \qquad : K_{\CM} \rightarrow K_{\CM \times S}  \label{eqn:e k-theory} \\
&f_n = (\pi_- \times \pi_S)_* \Big(L^{n-r} (-1)^r \det \CU \cdot \pi_+^* \Big) : K_{\CM} \rightarrow K_{\CM \times S}  \label{eqn:f k-theory}
\end{align}
for all $n \in \BZ$, where we write $L = [\CL]$ for the $K$-theory class of the tautological line bundle \eqref{eqn:tautological} \footnote{Moreover, in \eqref{eqn:f k-theory} we let $\det \CU$ denote the determinant of the universal sheaf on $\fZ$ which parameterizes either of $\CF_x$ and $\CF'_x$ in the notation \eqref{eqn:simple}. Since the coherent sheaves $\CF$ and $\CF'$ only differ in codimension 2, they have isomorphic determinant.}. Let $q=[\CK_S] \in K_{S}$, and consider the composition
\begin{equation}
\label{eqn:h k-theory series}
h^\pm(z) = \sum_{m=0}^{\infty} \frac {h^\pm_m}{z^{\pm m}} : K_{\CM} \xrightarrow{\text{pull-back}} K_{\CM \times S} \xrightarrow{\otimes \wedge^\bullet\left(\frac {z(q-1)}{\CU}\right)} K_{\CM \times S}
\end{equation}
where we expand the rational function $\wedge^\bullet(\frac {z(q-1)}{\CU})$ in powers of $z^{\mp 1}$. We will package the operators \eqref{eqn:e k-theory} and \eqref{eqn:f k-theory} into power series
\begin{align}
&e(z) = \sum_{n=-\infty}^{\infty} \frac {e_n}{z^{n}} = (\pi_+ \times \pi_S)_*\left[\delta \left(\frac {L}z \right) \cdot \pi_-^* \right] \label{eqn:e k-theory series} \\
&f(z) = \sum_{n=-\infty}^{\infty} \frac {f_n}{z^{n}} = (\pi_- \times \pi_S)_*\left[\delta \left(\frac {L}z \right) \frac {\det \CU}{(-z)^r} \cdot \pi_+^* \right] \label{eqn:f k-theory series}
\end{align}

\medskip

\subsection{}

Let us formally write
\begin{equation}
\label{eqn:chern roots k-theory}
[\Omega_S^1] = q_1+q_2 \quad \text{and} \quad [\Omega_S^2] = q = q_1q_2
\end{equation}
While $q_1$ and $q_2$ are not themselves elements of $K_S$, any symmetric polynomial in them is a well-defined element in $K_S$. We are now ready to give the commutation relations between the operators $e_n, f_n, h_m^\pm$. We will find that these relations look best when expressed in terms of the series $e(z), f(z), h^\pm(z)$. As in Subsection \ref{sub:compositions}, let us color the factors of $\textcolor{red}{S} \times \textcolor{blue}{S}$ in red and blue, and use the notation $\textcolor{red}{q}$ and $\textcolor{blue}{q}$ for the pull-back of the class $q$ from either one of the factors to $\textcolor{red}{S} \times \textcolor{blue}{S}$. Consider
\begin{equation}
\label{eqn:def zeta k-theory}
\zeta^\kth(x) = \wedge^\bullet(-x \CO_\Delta) = 1 + \frac {x [\CO_\Delta]}{(1-x)(1-xq)} \in K_{\textcolor{red}{S} \times \textcolor{blue}{S}}(x)
\end{equation}
where $\Delta \hookrightarrow S \times S$ is the class of the diagonal \footnote{For a proof of the second equality in \eqref{eqn:def zeta k-theory}, see \cite[Proposition 5.24]{N shuffle surf}.}, and 
\begin{align}
&\tzeta^{\kth}_+(x) = \zeta^{\kth}(x) \Big(1- x\textcolor{red}{q}\Big)\left(1- \frac {\textcolor{blue}{q}}x\right) \in \frac {K_{\textcolor{red}{S} \times \textcolor{blue}{S}}[x^{\pm 1}]}{1-x} \label{eqn:def zeta k-theory plus} 
 \\
&\tzeta^{\kth}_-(x) = \zeta^{\kth}(x) \left(1- \frac {\textcolor{red}{q}}x\right)\left(1- x\textcolor{blue}{q} \right) \in \frac {K_{\textcolor{red}{S} \times \textcolor{blue}{S}}[x^{\pm 1}]}{1-x} \label{eqn:def zeta k-theory minus} 
\end{align}
In \eqref{eqn:def zeta k-theory}, the symbol $q$ in the denominator denotes either $\textcolor{red}{q}$ or $\textcolor{blue}{q}$; it is not important which is chosen due to the presence of $[\CO_\Delta]$. In what follows, for any operators $x,y : K_{\CM} \rightarrow K_{\CM \times S}$, we will consider the following compositions of operators
\begin{align}
&\textcolor{red}{x} \textcolor{blue}{y} : K_{\CM} \xrightarrow{\textcolor{blue}{y}} K_{\CM \times \textcolor{blue}{S}} \xrightarrow{\textcolor{red}{x} \boxtimes \text{Id}_{\textcolor{blue}{S}}} K_{\CM\times \textcolor{red}{S} \times \textcolor{blue}{S}} \label{eqn:composition xy k-theory} \\
& \textcolor{blue}{y} \textcolor{red}{x} : K_{\CM} \xrightarrow{\textcolor{red}{x}} K_{\CM \times \textcolor{red}{S}} \xrightarrow{\textcolor{blue}{y} \boxtimes \text{Id}_{\textcolor{red}{S}}} K_{\CM\times \textcolor{red}{S} \times \textcolor{blue}{S}} \label{eqn:composition yx k-theory}
\end{align}
by analogy with the similar notions in Subsection \ref{sub:compositions}.

\medskip

\begin{theorem}
\label{thm:k-theory}

(\cite{N shuffle surf}) We have the following equalities of operators $K_{\CM} \rightarrow K_{\CM \times \textcolor{red}{S} \times \textcolor{blue}{S}}$
\begin{equation}
\label{eqn:rel 1 k-theory}
\textcolor{red}{e(z)} \textcolor{blue}{e(w)} \tzeta^{\ekth}_{-}\left(\frac wz\right) =  \textcolor{blue}{e(w)}\textcolor{red}{e(z)}  \tzeta^{\ekth}_{+}\left(\frac zw\right) 
\end{equation}
\begin{equation}
\label{eqn:rel 2 k-theory}
 \textcolor{blue}{f(w)} \textcolor{red}{f(z)} \tzeta^{\ekth}_{-}\left(\frac wz\right) = \textcolor{red}{f(z)} \textcolor{blue}{f(w)} \tzeta^{\ekth}_{+}\left(\frac zw\right) 
\end{equation}
\begin{equation}
\label{eqn:rel 3 k-theory}
\textcolor{red}{h^\pm(z)} \textcolor{blue}{e(w)} =  \textcolor{blue}{e(w)} \textcolor{red}{h^\pm(z)} \frac {\zeta^{\ekth}\left(\frac zw\right)}{\zeta^{\ekth}\left(\frac wz\right)}
\end{equation}
\begin{equation}
\label{eqn:rel 4 k-theory}
\textcolor{blue}{f(w)} \textcolor{red}{h^\pm(z)} =\textcolor{red}{h^\pm(z)} \textcolor{blue}{f(w)} \frac {\zeta^{\ekth}\left(\frac zw\right)}{\zeta^{\ekth}\left(\frac wz\right)}
\end{equation}
\begin{equation}
\label{eqn:rel 5 k-theory}
\left [\textcolor{red}{f(z)}, \textcolor{blue}{e(w)} \right] = \delta \left(\frac zw\right)  \cdot \Delta_* \left( \frac{h^+(z)-h^-(w)}{1-q}\right)
\end{equation}
as well as $[\textcolor{red}{h^\pm(z)}, \textcolor{blue}{h^{\pm'}(w)}] = 0$, for all $\pm,\pm' \in \{+,-\}$. The class $1-q$ in the denominator of \eqref{eqn:rel 5 k-theory} does not pose an issue because $h_0^+-h_0^-$ and all $\{h^\pm_m\}_{m>0}$ are multiples of $1-q$, as is clear from \eqref{eqn:h k-theory series}. 

\end{theorem}

\medskip

\subsection{}

To unravel the formulas in Theorem \ref{thm:k-theory}, we first note that \eqref{eqn:rel 5 k-theory} reads
\begin{equation}
\label{eqn:rel 5 k-theory explicit}
[\textcolor{red}{f_n}, \textcolor{blue}{e_m}] = \Delta_* \left( \frac {1}{1-q} \cdot \begin{cases} h^+_{m+n} &\text{if }m+n > 0 \\ h^+_0 - h^-_0 &\text{if }m+n =0 \\ -h^-_{-m-n}&\text{if }m+n < 0 \end{cases} \right)
\end{equation}
Let us now unravel relations \eqref{eqn:rel 3 k-theory} and \eqref{eqn:rel 4 k-theory}. We start by noting that 
$$
\frac {\zeta^{\kth}\left(\frac zw\right)}{\zeta^{\kth}\left(\frac wz\right)} = \frac {\displaystyle 1 + [\CO_\Delta] \cdot \frac {zw}{(z-w)(zq-w)}}{\displaystyle 1 + [\CO_\Delta] \cdot \frac {zw}{(z-w)(z-wq)}} = 1 + \sum_{a = 1}^{\infty} \Delta_*(\gamma_{a}^\pm) \frac {w^{\pm a}}{z^{\pm a}}
$$
where $\gamma_{a}^\pm \in K_S$ is the same polynomial expression in the classes $q_1+q_2$ and $q = q_1q_2$ of \eqref{eqn:chern roots k-theory} as the one we already encountered in \eqref{eqn:expansion toroidal}. These classes arise because of the identity $[\CO_\Delta]  [\CO_\Delta] = (1-q_1)(1-q_2)[\CO_\Delta]$ in $K_{\textcolor{red}{S} \times \textcolor{blue}{S}}$. Thus, relations \eqref{eqn:rel 3 k-theory} and \eqref{eqn:rel 4 k-theory} take the following form, for all $m \in \BZ$ and $n \geq 0$
\begin{equation}
\label{eqn:rel 3 k-theory explicit}
[\textcolor{red}{h^\pm_{n}}, \textcolor{blue}{e_{m}}] = \Delta_* \left(  \sum_{a = 1}^{n} \gamma_{a}^\pm  \underbrace{\textcolor{blue}{e_{m\pm a}} \textcolor{red}{h_{n-a}^\pm} \Big|_\Delta}_{\text{an operator } K_{\CM} \rightarrow K_{\CM \times \textcolor{red}{S} \times \textcolor{blue}{S}} \xrightarrow{|_\Delta} K_{\CM \times S}} \right)
\end{equation}
\begin{equation}
\label{eqn:rel 4 k-theory explicit}
[\textcolor{blue}{f_{m}}, \textcolor{red}{h^\pm_{n}}] = \Delta_* \left(  \sum_{a =1}^{n} \gamma^\pm_a  \underbrace{\textcolor{red}{h^\pm_{n-a}}\textcolor{blue}{f_{m\pm a}}  \Big|_\Delta}_{\text{an operator } K_{\CM} \rightarrow K_{\CM \times \textcolor{red}{S} \times \textcolor{blue}{S}} \xrightarrow{|_\Delta} K_{\CM \times S}} \right)
\end{equation}
Compare the formulas above with \eqref{eqn:rel 3 tor explicit}--\eqref{eqn:rel 4 tor explicit}. 

\medskip

\noindent As for formulas \eqref{eqn:rel 1 k-theory}--\eqref{eqn:rel 2 k-theory}, let us first observe that
\begin{align*}
&\tzeta^{\kth}_{+}\left(\frac zw\right) = \left(1-\frac {z\textcolor{red}{q}}w\right)\left(1-\frac {w\textcolor{blue}{q}}z \right) - [\CO_\Delta] \frac {z-wq}{z-w} \\
&\tzeta^{\kth}_{-}\left(\frac wz\right) = \left(1-\frac {z\textcolor{red}{q}}w \right)\left(1-\frac {w\textcolor{blue}{q}}z \right) - [\CO_\Delta] \frac {zq-w}{z-w} 
\end{align*}
Therefore, \eqref{eqn:rel 1 k-theory} and \eqref{eqn:rel 2 k-theory} are equivalent to the following formulas for all $m,n \in \BZ$
\begin{multline}
\label{eqn:rel 1 k-theory explicit}
[\textcolor{red}{e_{n+3}},\textcolor{blue}{e_m}] - \left(\frac 1{\textcolor{red}{q}} + 1 + \textcolor{blue}{q}\right) [\textcolor{red}{e_{n+2}},\textcolor{blue}{e_{m+1}}] + \\  \left(\frac 1{\textcolor{red}{q}} + \frac {\textcolor{blue}{q}}{\textcolor{red}{q}}+  \textcolor{blue}{q}\right) [\textcolor{red}{e_{n+1}},\textcolor{blue}{e_{m+2}}] -\frac {\textcolor{blue}{q}}{\textcolor{red}{q}} [\textcolor{red}{e_{n}},\textcolor{blue}{e_{m+3}}] +  
\end{multline}
$$
+ \Delta_*\left(e_{n+2} e_{m+1} \Big|_\Delta - \frac 1q e_{n+1} e_{m+2} \Big|_\Delta - \frac 1q e_{m+1} e_{n+2} \Big|_\Delta + e_{m+2}e_{n+1} \Big|_\Delta\right) = 0
$$
and 
\begin{multline}
\label{eqn:rel 2 k-theory explicit}
[\textcolor{blue}{f_m},\textcolor{red}{f_{n+3}}] - \left( \frac 1{\textcolor{red}{q}} + 1 + \textcolor{blue}{q}\right)  [\textcolor{blue}{f_{m+1}},\textcolor{red}{f_{n+2}}] + \\  \left(\frac 1{\textcolor{red}{q}} + \frac {\textcolor{blue}{q}}{\textcolor{red}{q}}+  \textcolor{blue}{q}\right) [\textcolor{blue}{f_{m+2}},\textcolor{red}{f_{n+1}}] - \frac {\textcolor{blue}{q}}{\textcolor{red}{q}} [\textcolor{blue}{f_{m+3}},\textcolor{red}{f_{n}}] +  
\end{multline}
$$
+ \Delta_*\left(f_{m+1} f_{n+2}  \Big|_\Delta - \frac 1q  f_{m+2} f_{n+1}\Big|_\Delta - \frac 1q  f_{n+2} f_{m+1} \Big|_\Delta + f_{n+1} f_{m+2}\Big|_\Delta\right) = 0
$$
The discussion in Subsections \ref{sub:diagonal 1} and \ref{sub:diagonal 2} applies verbatim to the situation of $K$-theory. In particular, we have the notion of operators $a,b : K_{\CM} \rightarrow K_{\CM \times S}$ having diagonal commutator, in which case we will write
$$
[\textcolor{red}{a}, \textcolor{blue}{b}] = \Delta_*\left([\textcolor{red}{a}, \textcolor{blue}{b}]_{\text{red}}\right)
$$
All the operators $e_n,f_n,h_m^\pm$ of \eqref{eqn:e k-theory}, \eqref{eqn:f k-theory}, \eqref{eqn:h k-theory series} pairwise have diagonal commutator, and the notion of an action analogous to Definition \ref{def:action}
\begin{equation}
\label{eqn:action toroidal}
\UU \curvearrowright K_{\CM}
\end{equation}
is well-defined (the only modification we need to make is that \eqref{eqn:commutator} must read
$$
[\Phi(x),\Phi(y)]_{\text{red}} = \Phi \left(\frac {[x,y]}{(1-q_1)(1-q_2)}\right)
$$
Note that the right-hand side is well-defined because of \eqref{eqn:divisible toroidal}). Then Theorem \ref{thm:k-theory} precisely states that the operators \eqref{eqn:e k-theory}, \eqref{eqn:f k-theory}, \eqref{eqn:h k-theory series} yield an action \eqref{eqn:action toroidal}, which is the content of Theorem \ref{thm:k-theory intro}.

\bigskip

\section{Elliptic cohomology}

\medskip

\subsection{}

In the present Section, we will recall the elliptic cohomology of smooth varieties over $\BC$ and elliptic quantum groups, and connect these two notions by proving analogues of Theorems \ref{thm:coh} and \ref{thm:k-theory}. Given a smooth algebraic variety $X$, recall that its cohomology $H_X$ is a ring endowed with a Chern polynomial
\begin{equation}
\label{eqn:chern coh}
c(\CV,z) = (z-\ell_1)\dots(z-\ell_r) \in H_X[z]
\end{equation}
for any rank $r$ locally free sheaf $\CV$ on $X$. The symbols $\ell_1,\dots,\ell_r$ are called the Chern roots of $\CV$, and while they are not individually well-defined, any symmetric polynomial in $\ell_1,\dots,\ell_r$ is a well-defined element of $H_X$. In other words, specifying the Chern polynomial of $\CV$ is tantamount to specifying the ring homomorphism
$$
\BC[x_1,\dots,x_r]^{\text{Sym}} \longrightarrow H_X
$$
determined by $x_i \mapsto \ell_i$. Dually, this amounts to a map of schemes
$$
\spec(H_X) \rightarrow (\BA^1)^{(r)} 
$$
Above, $C^{(r)}$ denotes the $r$-th symmetric power of any algebraic curve $C$, whose closed points are unordered $r$-tuples of closed points of $C$. The situation of $K$-theory is analogous to that of cohomology, except that we replace $\BA^1$ by $\BA^1 \backslash \{0\}$.

\medskip

\begin{definition}
\label{def:gkv}

(\cite{gkv}) Let $E$ be an elliptic curve over $\BC$. The (0-th) elliptic cohomology associated to $E$ is a contravariant functor
\begin{equation}
\label{eqn:def elliptic}
\Big( \emph{smooth varieties over }\BC \Big) \quad \xrightarrow{X \mapsto \eellcoh_X} \quad \emph{Commutative rings}
\end{equation}
The (total) Chern class of a rank $r$ locally free sheaf $\CV$ is a map of schemes
\begin{equation}
\label{eqn:chern elliptic}
\espec(\eellcoh_X) \xrightarrow{c_{\CV}} E^{(r)}
\end{equation}
satisfying the following properties with respect to direct sums and tensor products
$$
\xymatrixcolsep{3pc}\xymatrix{
\espec(\eellcoh_X) \ar[dr]_-{c_{\CV \oplus \CV'}} \ar[r]^-{c_{\CV} \times c_{\CV'}} & E^{(r)} \times E^{(r')} \ar[d]^-{\oplus} \\
 &E^{(r+r')}} \qquad \xymatrixcolsep{3pc}\xymatrix{
\espec(\eellcoh_X) \ar[dr]_-{c_{\CV \otimes \CV'}} \ar[r]^-{c_{\CV} \times c_{\CV'}} & E^{(r)} \times E^{(r')} \ar[d]^-{\otimes} \\
 &E^{(rr')}}
$$
for any locally free sheaves $\CV, \CV'$ of ranks $r,r'$, respectively. In the diagrams above, the vertical arrows are the maps given in terms of closed points by
\begin{equation}
\label{eqn:oplus}
(x_1,\dots,x_r) \oplus (y_1,\dots,y_{r'}) = (x_1,\dots,x_r, y_1,\dots,y_{r'})
\end{equation}
\begin{equation}
\label{eqn:otimes}
(x_1,\dots,x_r) \otimes (y_1,\dots,y_{r'}) = (\dots, x_i y_j, \dots)_{1\leq i \leq r, 1 \leq j \leq r'}
\end{equation}

\end{definition}

\medskip

\noindent The reason why $E$ has to be an elliptic curve is the presence of the group law of $E$ in the formula for the map $\otimes$ in \eqref{eqn:otimes}. It is also why $\BA^1$ (with the additive group law) and $\BA^1 \backslash 0$ (with the multiplicative group law) arise in the situations of cohomology and $K$-theory, respectively.

\medskip

\subsection{}

As explained in \cite{gkv}, the fact that the target of \eqref{eqn:chern elliptic} is not affine is one of the reasons one would use the scheme
$$
A_X = \spec(\ellcoh_X)
$$
in order to encode elliptic cohomology (this is all the more because of the analogous situation of $T$-equivariant cohomology for a torus $T = (\BC^*)^k$, in which case the analogue of $\ellcoh_X$ is a sheaf of commutative rings over $E^k$, and thus the analogue of $A_X$ is an affine-over-projective scheme). The contravariance of the functor \eqref{eqn:def elliptic} entails the existence of a ring homomorphism
$$
f^* : \ellcoh_Y \rightarrow \ellcoh_X
$$
or dually, a map of schemes
$$
\widetilde{f} : A_X \rightarrow A_Y
$$
for any morphism $f : X \rightarrow Y$. As a consequence of the axioms developed in \cite{gkv}, for any rank $r$ vector bundle $\CV$ on $X$, the elliptic cohomology of the projective bundle
$$
\BP_X(\CV) = \text{Proj}_X(\text{Sym}^\bullet\CV) \xrightarrow{\pi} X
$$
has the property that the following square is Cartesian
\begin{equation}
\label{eqn:cartesian}
\xymatrix{
A_{\BP_X(\CV)} \ar[d]_{\widetilde{\pi}} \ar[r] & E \times E^{(r-1)} \ar[d]^{\text{the map \eqref{eqn:oplus}}} \\
A_X \ar[r]^{c_{\CV}} & E^{(r)}}
\end{equation}
Composing the top arrow above with projection to the first factor yields a map $A_{\BP_X(\CV)} \rightarrow E$, which is none other than $c_{\CO(1)}$, i.e. the map \eqref{eqn:chern elliptic} applied to the tautological line bundle $\CO(1)$ on $\BP_X(V)$. Formula \eqref{eqn:cartesian} is the scheme version of the well-known isomorphism in ordinary cohomology
$$
H_{\BP_X(\CV)} \cong H_X[z] / c(\CV,z)
$$

\medskip

\subsection{} 
\label{sub:theta functions}

Ordinary cohomology classes are elements of $H_X$, or equivalently, functions on $\spec(H_X)$. Meanwhile, elliptic cohomology classes are more naturally construed as sections of line bundles on $A_X$. The prime example of this comes from the map \eqref{eqn:chern elliptic}: there are no non-constant functions on $E^{(r)}$, but there are plenty of meromorphic functions (i.e. sections of line bundles) on $E^{(r)}$. Pulling back the aforementioned meromorphic functions via $c_{\CV}$ naturally yields sections of line bundles on $A_X$, or in other words, elements of locally free rank 1 modules over $\ellcoh_X$.

\medskip

\begin{remark} If one works over $\BC$, the ring $\eellcoh_X$ is Artinian for any smooth variety $X$, so any locally free rank 1 module is isomorphic to $\eellcoh_X$. However, there is no canonical choice of such an isomorphism; this means that while any elliptic cohomology class can be construed as an element in $\eellcoh_X$ ``up to a unit", there is no canonical choice for this unit. Therefore, it is still beneficial to think of elliptic cohomology classes as sections of line bundles on $A_X$, and not functions on $A_X$.

\end{remark}

\medskip

\noindent We will now explicitly recall theta functions, which are sections of line bundles on the elliptic curve $E$. To do this, let us consider the multiplicative presentation
$$
E = \BC^* / p^{\BZ}
$$
where $p$ is a complex number with $0 < |p| < 1$ (in the more common presentation $E = \BC/(\BZ \oplus \BZ \tau)$, we have $p = e^{2\pi i \tau}$). Then for any pair $(n,\lambda) \in \BZ \times \BC^*$, which we will refer to as factors of automorphy, we may define the line bundle
$$
\DD_{n,\lambda} = \Big\{\text{meromorphic functions }f(z) \text{ on }\BC^*, \text{ such that } \frac {f(z)}{f(zp)} = z^n \lambda, \forall z \in \BC^* \Big\}
$$
on $E$. It is easy to see that multiplication of functions induces an isomorphism
$$
\DD_{n,\lambda} \otimes \DD_{m,\mu} \xrightarrow{\sim} \DD_{n+m,\lambda\mu}
$$
for all $n,m \in \BZ$ and $\lambda, \mu \in \BC^*$, and that $\DD_{0,1}$ is the trivial line bundle. In particular, for all $k \in \BZ$, multiplication by $z^k$ induces an isomorphism
$$
\DD_{n,\lambda} \xrightarrow{\sim} \DD_{n,\lambda p^{-k}}
$$
We will interpret the expressions $z^k$ as sections of the different (albeit isomorphic) line bundles $\DD_{0,p^{-k}}$ on $E$. We will refer to $z$ as the standard coordinate on $E$, although as explained in the preceding sentence, this is an abuse of terminology.

\medskip

\begin{proposition}
\label{prop:theta}

The Jacobi theta function
\begin{equation}
\label{eqn:def theta}
\vartheta(z) = (1-z) \prod_{s=1}^{\infty} \frac {(1-p^s z)(1-p^s z^{-1})}{(1-p^s)^2}
\end{equation}
is a section of $\DD_{1,-1}$.

\end{proposition}

\medskip

\noindent The non-standard choice of normalization in \eqref{eqn:def theta} is due to the fact that it ensures
\begin{equation}
\label{eqn:zero}
\vartheta(z) \Big|_{p=0} = 1-z 
\end{equation}
and
\begin{equation}
\label{eqn:residue} 
\underset{z=1}{\text{Res}} \ \frac 1{\vartheta(z)} = -1
\end{equation}
Proposition \ref{prop:theta} above is an immediate consequence of the simple identity
\begin{equation}
\label{eqn:theta quasiper}
\vartheta(zp) = -\frac {\vartheta(z)}z \qquad \Leftrightarrow \qquad \vartheta \left(\frac zp\right) = - \frac {z \vartheta(z)}p 
\end{equation}
More generally, for any $x_1,\dots,x_n,y_1,\dots,y_m \in \BC^*$, the ratio
\begin{equation}
\label{eqn:ratio of thetas}
\frac {\displaystyle \vartheta\left(zx_1\right) \dots \vartheta\left(zx_n\right)}{\displaystyle \vartheta\left(zy_1\right) \dots \vartheta\left(zy_m\right)}
\end{equation}
is a global section of $\DD_{n-m,(-1)^{n-m} \frac {x_1 \dots x_n}{y_1\dots y_m}}$. Finally, we point out the formula
\begin{equation}
\label{eqn:change}
\vartheta\left(z^{-1}\right) = -z^{-1} \vartheta(z)
\end{equation}

\medskip

\subsection{} 
\label{sub:push 1}

Non-trivial line bundles (which we will refer to as ``twists") arise in elliptic cohomology when defining push-forwards. To see this, if $\iota : X \hookrightarrow Y$ is the zero locus of a regular section of a rank $r$ locally free sheaf $\CV$ on $Y$, then we have
$$
\iota_*(1) = \ell_1 \dots \ell_r \in H_Y
$$
in usual cohomology, where $\ell_1,\dots,\ell_r$ are the Chern roots of $\CV$. Similarly, we have 
$$
\iota_*(1) = \left(1-L_1^{-1} \right)\dots \left(1-L_r^{-1} \right) \in K_Y
$$
in $K$-theory, where we formally decompose $[\CV]$ as a sum of line bundles $L_1+\dots+L_r$ \footnote{This is simply an artifice; the rigorous definition of the $L_i$'s is such that $$\wedge^\bullet\left(\frac {\CV}z\right) = \left(1-\frac {L_1}z\right) \dots  \left(1-\frac {L_r}z\right)$$}. Generalizing the formula above, we would like to have in elliptic cohomology
\begin{equation}
\label{eqn:push elliptic}
`` \iota_*(1) = \vartheta \left(L^{-1}_1\right)\dots \vartheta \left(L^{-1}_r \right) \in \ellcoh_Y "
\end{equation}
However, $L_1,\dots,L_r$ here should be interpreted as the pull-backs under $c_{\CV}$ of the standard coordinates on $E^{(r)}$. Therefore, the right-hand side of \eqref{eqn:push elliptic} should be interpreted not as a function on $A_Y = \text{Spec}(\text{Ell}_Y)$, but as the section
$$
\vartheta\left(\CV^\vee \right) := \vartheta \left(L^{-1}_1\right)\dots \vartheta\left(L^{-1}_r\right)
$$
of the following line bundle on $A_Y$
$$
\Theta \left(\CV \right) := c_{\CV^\vee}^* \left( (\DD_{1,-1} \boxtimes \dots \boxtimes \DD_{1,-1})^{\text{sym}}\right)
$$
(in the right-hand side, we descend the tensor product of the line bundles of Proposition \ref{prop:theta} from the $r$ factors of $E^r$ to a line bundle on $E^{(r)}$). Note that our $\Theta(\CV)$ is actually $\Theta(-\CV^\vee)$ in the notation of \cite{gkv}. 

\medskip

\begin{proposition}

For any $\CV,\CW$, there is an isomorphism of line bundles 
$$
\Theta(\CV\oplus \CW) \cong \Theta(\CV) \otimes \Theta(\CW)
$$ 
Therefore, there exists a group homomorphism
$$
K_Y \xrightarrow{\Theta} \emph{Pic}(A_Y)
$$
\begin{equation}
\label{eqn:theta k}
\Theta(\CV-\CW) = \Theta(\CV) \otimes \Theta(\CW)^{-1}
\end{equation}
Moreover, the line bundle \eqref{eqn:theta k} has the section
$$
\vartheta(\CV^\vee-\CW^\vee) = \frac {\vartheta(L^{-1}_1) \dots \vartheta(L^{-1}_r)}{\vartheta(M^{-1}_1) \dots \vartheta(M^{-1}_s)}
$$
for any vector bundles $[\CV] = L_1 + \dots + L_r$ and $[\CW] = M_1 + \dots + M_s$ on $Y$.

\end{proposition}

\medskip

\noindent Thus, formula \eqref{eqn:push elliptic} should be read as
\begin{equation}
\label{eqn:push closed embedding}
\iota_*(1) = \vartheta\left( \CV^\vee \right) \text{ as a section of } \Theta \left( \CV \right)
\end{equation}

\medskip

\subsection{} 
\label{sub:push 2}

Motivated by the discussion above, the general construction is the following. For a smooth algebraic variety $X$, we will denote its tangent bundle by $\text{Tan}_X$.

\medskip

\begin{definition}
\label{def:push}

(\cite{gkv}) For a proper map $f : X \rightarrow Y$, push-forward is a morphism 
\begin{equation}
\label{eqn:def push}
f_* : \widetilde{f}_*(\Theta(\emph{Tan}_X - f^*(\emph{Tan}_Y)) \rightarrow \CO
\end{equation}
of coherent sheaves on $A_Y$ (the latter condition encodes the projection formula).

\end{definition}

\medskip

\noindent Since all the elliptic cohomology schemes $A_X$ in the present paper will be affine, we simplify \eqref{eqn:def push} by always thinking about push-forward as a map of $\ellcoh_Y$-modules
$$
f_* : \Big(\text{a rank 1 locally free }\ellcoh_X \text{-module} \Big) \rightarrow \ellcoh_Y
$$
Moreover, in all our computations we will abuse notation and write
\begin{equation}
\label{eqn:abuse}
f_* : \ellcoh_X \rightarrow \ellcoh_Y
\end{equation}
because the formula for $f_*$ will always be given by a product of $\vartheta$ functions, which completely encodes which rank 1 locally free $\ellcoh_X$-module is the correct domain of $f_*$. For example, if $\iota : X \hookrightarrow Y$ is the zero locus of a regular section of a vector bundle $\CV$ on $Y$, then the push-forward is explicitly given by
\begin{equation}
\label{eqn:push embedding}
\iota_*(\widetilde{\sigma}) = \sigma \cdot \vartheta(\CV^\vee) 
\end{equation}
for any element $\sigma$ in a rank 1 locally free $\ellcoh_Y$-module, where
$$
\widetilde{\sigma} = \widetilde{\iota}^*(\sigma)
$$
denotes the pull-back of $\sigma$ to a rank 1 locally free $\ellcoh_X$-module. Formula \eqref{eqn:push embedding} encodes the fact that the domain and codomain of $\iota_*$ are rank 1 locally free modules which differ by a twist with $\Theta(-\iota^*(\CV)) = \Theta(\text{Tan}_X - \iota^*(\text{Tan}_Y))$. 

\medskip

\subsection{} Another case of great interest to us is push-forward along a projective bundle
$$
\xymatrix{
\BP_X(\CV) \ar[d]^{\pi} \\
X}
$$
We have the map of schemes induced by the Cartesian diagram \eqref{eqn:cartesian} 
$$
A_{\BP_X(V)} \xrightarrow{\widetilde{\pi} \times c_{\CO(1)}} A_X \times E
$$
In this case, for any section $\sigma(z)$ of a line bundle on $A_X \times E$ (we write $z$ for the coordinate on the latter copy of $E$), the following formula was proved in \cite{zz}:
\begin{equation}
\label{eqn:push projective bundle}
\pi_*(\widetilde{\sigma}(\taut)) = \sum_{i=1}^r \frac {\sigma(L_i)}{\prod_{j \neq i} \vartheta\left(\frac {L_j}{L_i}\right)}
\end{equation}
where $L_1,\dots,L_r$ denote the pullbacks under $c_{\CV}$ of the coordinates on $E^{(r)}$ and
$$
\widetilde{\sigma}(\taut) = (\widetilde{\pi} \times c_{\CO(1)})^*(\sigma(z))
$$
(the notation in the left-hand side is motivated by the fact that $\taut := [\CO(1)]$ is the pull-back of the coordinate $z$ under the map $c_{\CO(1)} : A_{\BP_X(V)} \rightarrow E$). Formula \eqref{eqn:push projective bundle} encodes the fact that the domain and codomain of $\pi_*$ are sections of line bundles which differ by a twist with
$$
\Theta \left(\pi^*(\CV^\vee) \otimes \CO_{\BP_X(\CV)}(1) - \CO_{\BP_X(\CV)}\right)\ = \Theta\left(\text{Tan}_{\BP_X(\CV)} - \pi^*(\text{Tan}_X) \right) 
$$
For the purpose of computations, we will assume that $L_1,\dots,L_r$ in formula \eqref{eqn:push projective bundle} are complex numbers of absolute value in $(|p| ,1]$. This is not really too far from the truth, since they are the pull-backs to $\ellcoh_X$ of the usual coordinates on $E = \BC^*/p^{\BZ}$, and the latter can always be realized in the chosen interval up to multiplication by some power of $p$. Then because the function $\vartheta(z)$ has a unique simple pole at $z=1$ with residue $-1$, formula \eqref{eqn:push projective bundle} can be rewritten as follows
\begin{equation}
\label{eqn:push explicit 1}
\pi_*(\widetilde{\sigma}(\taut)) = \int_{1-p} \frac {\sigma(z)}{\displaystyle \vartheta\left(\frac {\CV}z \right)} 
\end{equation}
where in the RHS we set 
$$
\vartheta\left(\frac {\CV}z \right) = \prod_{i=1}^r \vartheta \left(\frac {L_i}z \right)
$$
and for any meromorphic function $F(z)$, we write
\begin{equation}
\label{eqn:def integral}
\int_{1-p} F(z) = \int_{|z|=1} F(z) \frac {dz}{2\pi i z} - \int_{|z|=|p|} F(z) \frac {dz}{2\pi i z}
\end{equation}

\medskip

\subsection{} 
\label{sub:push 3}

We will also be interested in a slightly more general kind of projectivization, which combines the situations of closed embeddings and projective bundles that were discussed above. Explicitly, suppose we have a short exact sequence on $X$
$$
0 \rightarrow \CW \rightarrow \CV \rightarrow \CE \rightarrow 0
$$
where $\CV$ and $\CW$ are locally free sheaves, but $\CE$ is simply a coherent sheaf. Then we have a commutative diagram
\begin{equation}
\label{eqn:proj sheaf}
\xymatrix{
\BP_X(\CE) \ar@{^{(}->}[r]^\iota \ar[dr]_{\pi'} & \BP_X(\CV) \ar[d]^{\pi} \\
& X}
\end{equation}
with the closed embedding $\iota$ defined as the zero locus of the following composition
$$
\pi^*(\CW) \rightarrow \pi^*(\CV) \rightarrow \CO(1)
$$
We assume that the composition above is a regular section of $\pi^*(\CW^\vee) \otimes \CO(1)$, which is equivalent to the condition that $\BP_X(\CE)$ has expected dimension equal to $\text{rank }\CE - 1$ over $X$. In this case, we may compute $\pi'_*$ by combining formulas \eqref{eqn:push embedding} and \eqref{eqn:push explicit 1}
\begin{equation}
\label{eqn:push explicit 2}
\begin{split} \pi'_*(\widetilde{\sigma}(\taut)) &=  \int_{1-p} \frac {\displaystyle \sigma(z)\vartheta\left(\frac {\CW}z\right)}{\displaystyle \vartheta\left(\frac {\CV}z \right)} \\ &= \int_{1-p} \frac {\sigma(z)}{\displaystyle \vartheta\left(\frac {\CE}z\right)}
\end{split}
\end{equation}

\medskip

\begin{remark}
\label{rem:fac auto}

A particularly important case of \eqref{eqn:push explicit 2} is when $\sigma(z) = z^m$ for some integer $m$. We may assemble these special cases together by considering the series
$$
\delta \left(\frac zw \right) = \sum_{m = -\infty}^{\infty} \frac {z^m}{w^m}
$$
for some formal variable $w$. Then \eqref{eqn:push explicit 2} reads 
\begin{equation}
\label{eqn:push delta}
\pi'_* \left[ \displaystyle \delta \left(\frac {\emph{taut}}w \right) \right] = \frac 1{\displaystyle \vartheta\left(\frac {\CE}w \right)} \Big|_{1-p} 
\end{equation}
where for any meromorphic function $F(z)$, we write
\begin{equation}
\label{eqn:laurent expansion}
F(z) \Big|_{1-p} = \Big( \text{Laurent series of }F(z) \Big) - \Big( \text{Laurent series of }F(zp) \Big)
\end{equation}
All our Laurent series will be expanded near the circle $|z|=1$. Note that if $F(z)$ has poles between the circles $|z|=1$ and $|z|=|p|$, then the two Laurent series expansions in the right-hand side of \eqref{eqn:laurent expansion} can be quite different from each other.

\end{remark}

\medskip

\subsection{}

Let us now consider the moduli spaces $\CM$ of stable sheaves on a smooth projective surface $S$, with the notation in Subsection \ref{sub:basic moduli} and subject to Assumptions A and S. We will define elliptic cohomology versions of the operators in cohomology/$K$-theory from Subsections \ref{sub:operators cohomology} and \ref{sub:operators k-theory}, respectively. Recall the space $\fZ$ and the maps $\pi_\pm$ and $\pi_S$ from diagram \eqref{eqn:simple diag}, and let us define for all $n \in \BZ$
\begin{align}
&e_n = (\pi_+ \times \pi_S)_* \Big(L^n \cdot \pi_-^*\Big) \qquad \qquad \qquad : \ellcoh_{\CM} \rightarrow \ellcoh_{\CM \times S}  \label{eqn:e elliptic} \\
&f_n = (\pi_- \times \pi_S)_*\Big(L^{n-r} (-1)^r \det \CU \cdot \pi_+^* \Big) : \ellcoh_{\CM} \rightarrow \ellcoh_{\CM \times S} \label{eqn:f elliptic}
\end{align}
where $L$, $\det \CU$ denote the pull-backs under $c_{\CL}, c_{\det \CU} : \spec(\ellcoh_{\fZ}) \rightarrow E$ of the standard coordinate on $E$. As already explained, we are abusing notation in the formulas above, in that the domains and codomains of the operators $e_n,f_n$ are actually locally free rank 1 modules over the $\ellcoh$ rings in question, and not the rings themselves. However, it is straightforward to determine the specific twist by which these locally free rank 1 modules differ from each other by looking at the push-forwards $\pi_\pm \times \pi_S$, which were described in Proposition \ref{prop:projective bundles}; we will also give explicit formulas for the operators $e_n,f_n$ in terms of $\vartheta$ functions, from which the twist will be easy to read off. As before, it is helpful to package the operators \eqref{eqn:e elliptic} and \eqref{eqn:f elliptic} into power series
\begin{align}
&e(z) = \sum_{n=-\infty}^{\infty} \frac {e_n}{z^{n}} = (\pi_+ \times \pi_S)_*\left[\delta \left(\frac {L}z \right) \cdot \pi_-^* \right] \label{eqn:e elliptic series} \\
&f(z) = \sum_{n=-\infty}^{\infty} \frac {f_n}{z^{n}} = (\pi_- \times \pi_S)_*\left[\delta \left(\frac {L}z \right) \frac {\det \CU}{(-z)^r}\cdot \pi_+^* \right] \label{eqn:f elliptic series}
\end{align}
Finally, we define the following Laurent series (always expanded near $|z|=1$)
\begin{equation}
\label{eqn:h elliptic series}
\begin{split}
h^+(z) &= \sum_{n=-\infty}^{\infty} \frac {h^+_n}{z^n} : \ellcoh_{\CM} \xrightarrow{\text{pull-back}} \ellcoh_{\CM \times S} \xrightarrow{\cdot  \vartheta\left(\frac {z(q-1)}{\CU}\right)} \ellcoh_{\CM \times S} \\
h^-(z) &= \sum_{n=-\infty}^{\infty} \frac {h^-_n}{z^n} : \ellcoh_{\CM} \xrightarrow{\text{pull-back}} \ellcoh_{\CM \times S} \xrightarrow{\cdot \vartheta \left(\frac {zp(q-1)}{\CU}\right)} \ellcoh_{\CM \times S}
\end{split}
\end{equation}
\footnote{Note that $h_n^+ \neq p^n h_n^-$, because the Laurent series of a meromorphic function near the circle $|z|=1$ may be very different from the Laurent series of the same function near the circle $|z|=|p|$, due to the presence of poles between these two circles.} where $q = [\CK_S]$. Let us be more precise about the target of the series $h^\pm(z)$ above. Because $\CU^\vee(q-1)$ is a $K$-theory class of rank 0 and determinant equal to $q^r$, the expression
$$
\vartheta\left(\frac {z(q-1)}{\CU}\right)
$$
is a section of the pull-back of $\DD_{0,p^r}$ under the composition
$$
A_{\CM \times S} \xrightarrow{\widetilde{\text{proj}}} A_{S} \xrightarrow{c_{\CK_S}} E
$$

\medskip

\subsection{}

Having explained how to explicitly think of the operators $h_n^\pm$, let us do the same for the operators $e_n$ and $f_n$, and in the process elucidate the correct twist between their domain and codomain. In order to do so, we will recall universal classes. For any $k,l \in \BN$, let $\pi : \CM \times S^{k+l} \rightarrow \CM \times S^k$ and $\rho : \CM \times S^{k+l} \rightarrow S^{k+l}$ denote the standard projections and consider the Laurent series
\begin{equation}
\label{eqn:universal series}
\pi_* \left[ \vartheta \left(\frac {\CU_1}{z_1} \right) \dots \vartheta \left(\frac {\CU_{k+l}}{z_{k+l}} \right) \cdot \rho^*(\Gamma) \right] \in \ellcoh_{\CM \times S^k} 
\end{equation}
where $\CU_i$ denotes the universal sheaf \eqref{eqn:universal} pulled back from $\CM$ and the $i$-th factor of $S^{k+l}$, and $\Gamma \in \ellcoh_{S^{k+l}}$ is arbitrary. Any series coefficient of the expression above in the variables $z_1,\dots,z_{k+l}$ will be called a universal class, and will be denoted by
\begin{equation}
\label{eqn:universal class}
\Psi(\CU) \in \ellcoh_{\CM \times S^k} 
\end{equation}
We point out the abuse of notation in the formula above: the expression $\Psi(\CU)$ is not a function on $A_{\CM \times S^k}$, but more precisely a section of a line bundle on $A_{\CM \times S^k}$, determined by the particular choice of $\Gamma$ and the particular series coefficient in $z_1,\dots,z_{k+l}$ that we extract from \eqref{eqn:universal series}. 

\medskip

\begin{proposition}
\label{prop:e and f universal}

For any universal class $\Psi(\CU)$, we have the following formulas
\begin{equation}
\label{eqn:formula e universal}
e(z) \cdot \Psi(\CU) = \Psi(\CU + z \CO_\Delta) \vartheta \left(\frac {zq}{\CU} \right) \Big|_{1-p} 
\end{equation}
\begin{equation}
\label{eqn:formula f universal}
f(z) \cdot \Psi(\CU) = \Psi(\CU - z \CO_\Delta) \vartheta \left(- \frac z{\CU} \right) \Big|_{1-p} 
\end{equation}
(see the notation \eqref{eqn:laurent expansion} for $|_{1-p}$) where the diagonal $\Delta \subset S \times S$ identifies the factor of $S$ where the operators $e(z)$, $f(z)$ take values with the factor of $S$ where the universal sheaf $\CU$ is defined. Formulas for the operators $e_n$ and $f_n$ may be obtained by extracting the coefficient of $z^{-n}$ in the Laurent series in the right-hand sides.

\end{proposition}

\medskip

\noindent In the right-hand side of \eqref{eqn:formula e universal} and \eqref{eqn:formula f universal}, the notation
$$
\Psi\left(\CU + \beta \right)
$$
refers to replacing $\CU$ with $\CU + \beta$ in \eqref{eqn:universal series}, for any $\beta \in K_{\CM \times S^k}$.

\medskip

\begin{proof} Let us start by proving formula \eqref{eqn:formula e universal}. We have
\begin{equation}
\label{eqn:computation 1}
e(z) \cdot \Psi(\CU) = (\pi_+ \times \pi_S)_* \left( \delta \left(\frac Lz\right) \Psi(\CU) \right) 
\end{equation}
where in the right-hand side, $\Psi(\CU)$ denotes the universal class on $\fZ$ defined using the universal sheaf $\CU$ akin to \eqref{eqn:universal class}. However, on $\fZ \times S$ we actually have two universal sheaves $\CU$ and $\CU'$, due to the short exact sequence
$$
\xymatrix{0 \rightarrow \CU' \rightarrow \CU \rightarrow \CL \otimes \CO_{\Delta} \rightarrow 0 \ar@{.>}[d] \\ \fZ \times S}
$$
induced from the definition \eqref{eqn:simple}. Above, $\CU$ and $\CU'$ denote the universal sheaves parameterizing the stable sheaves $\CF$ and $\CF'$ from \eqref{eqn:simple} (respectively) and $\CO_\Delta$ denotes the pull-back of the diagonal under the homomorphism $\pi_S \times \text{Id}_S : \fZ \times S \rightarrow S \times S$.  Therefore, formula \eqref{eqn:computation 1} yields
\begin{equation}
\label{eqn:computation 2}
e(z) \cdot \Psi(\CU) = (\pi_+ \times \pi_S)_* \left[ \delta \left(\frac Lz\right)  \Psi(\CU' + \CL \otimes \CO_\Delta) \right]
\end{equation}
The following equality holds for any analytic function $f(z)$
$$
\delta\left(\frac Lz \right) f(L) = \delta\left(\frac Lz \right) f(z)
$$
and it allows us to rewrite \eqref{eqn:computation 2} as
\begin{equation}
\label{eqn:computation 3}
e(z) \cdot \Psi(\CU) = (\pi_+ \times \pi_S)_* \left[ \delta \left(\frac Lz\right)  \Psi(\CU' + z \CO_\Delta) \right]
\end{equation}
Meanwhile, as explained in the paragraph after Proposition \ref{prop:projective bundles}, the line bundle $L$ is the inverse of the tautological line bundle on 
$$
\fZ \cong \BP_{\CM' \times S}({\CU'}^\vee[1] \otimes \CK_S)
$$
Therefore, formula \eqref{eqn:push delta} may be used to evaluate \eqref{eqn:computation 3}, thus yielding \eqref{eqn:formula e universal} on the nose. As for \eqref{eqn:formula f universal}, the analogous treatment produces the following formula
\begin{equation}
\label{eqn:computation 4}
f(z) \cdot \Psi(\CU')  =  \Psi(\CU - z \CO_\Delta) \frac {\frac {\det \CU}{(-z)^r}}{\displaystyle \vartheta \left(\frac {\CU}z \right)} \Big|_{1-p}
\end{equation}
The fact that \eqref{eqn:computation 4} gives rise to formula \eqref{eqn:formula f universal} follows from the identity
$$
\vartheta \left(-\frac z{\CU} \right) = \frac {\frac {\det \CU}{(-z)^r}}{\displaystyle \vartheta \left(\frac {\CU}z \right)}
$$
which in turn is an immediate consequence of \eqref{eqn:change}.

\end{proof}

\medskip

\subsection{} We will now calculate the relations between the operators \eqref{eqn:e elliptic series}, \eqref{eqn:f elliptic series} and \eqref{eqn:h elliptic series}. Recall our convention of coloring the factors of $\textcolor{red}{S} \times \textcolor{blue}{S}$, as well as all elliptic cohomology classes on $\textcolor{red}{S} \times \textcolor{blue}{S}$ that are pulled back from one of the factors, in red and blue. Consider the elliptic analogue of the rational function \eqref{eqn:def zeta k-theory}
\begin{equation}
\label{eqn:def zeta elliptic}
\zeta^\ellcoh(x) = \vartheta(-x \CO_\Delta) \in \ellcoh_{\textcolor{red}{S} \times \textcolor{blue}{S}} (x)
\end{equation}
where $\CO_\Delta$ denotes the structure sheaf of the diagonal in $S \times S$. Because this structure sheaf has rank 0 and determinant 1, the coefficients of $\zeta^{\ellcoh}(x)$ actually lie in the ring $\ellcoh_{S \times S}$, and not in some locally free rank 1 module. Similarly, we have the following analogues of \eqref{eqn:def zeta k-theory plus} and \eqref{eqn:def zeta k-theory minus}
\begin{align}
&\tzeta^{\ellcoh}_+(x) = \zeta^{\ellcoh}(x) \vartheta \left(x \textcolor{red}{q} \right) \vartheta \left( \frac {\textcolor{blue}{q}}x \right)  \label{eqn:def zeta elliptic plus} 
 \\
&\tzeta^{\ellcoh}_-(x) = \zeta^{\ellcoh}(x) \vartheta \left(\frac {\textcolor{red}{q}}x \right) \vartheta \left( x \textcolor{blue}{q} \right)  \label{eqn:def zeta elliptic minus} 
\end{align}

\medskip

\begin{lemma}
\label{lem:integral}

The function $\zeta^{\eellcoh}(x)$ is of the form
\begin{equation}
\label{eqn:zeta want}
\zeta^{\eellcoh}(x) = 1 + \Delta_* \left( \frac {\text{quasi-periodic holomorphic function of }x}{\vartheta(x)\vartheta(xq)} \right)
\end{equation}
and so its only poles are at $x \in p^{\BZ}$ and $x  \in p^{\BZ}q^{-1}$. We have
\begin{equation}
\label{eqn:res zeta}
\underset{x=1}{\emph{Res}} \ \zeta^{\eellcoh}(x) = \Delta_*\left(\frac {-1}{\vartheta(q)} \right)
\end{equation}

\end{lemma}

\medskip

\noindent As a consequence of Lemma \ref{lem:integral}, the only poles of $\tzeta^{\ellcoh}_{\pm}(x)$ are $x \in p^{\BZ}$.

\medskip

\begin{proof} Because $\CO_\Delta$ is supported on the diagonal, we have $\zeta^{\ellcoh}(x)|_{S \times S \backslash \Delta} = 1$. Therefore, the excision long exact sequence in elliptic cohomology (\cite[(1.5.4)]{gkv}) implies that
$$
\zeta^{\ellcoh}(x) = 1 + \Delta_*(F(x))
$$
for some meromorphic function $F$ with coefficients in $\ellcoh_S$. It remains to show that $F(x)$ only has poles at $x \in p^{\BZ}$ and $x \in p^{\BZ} q^{-1}$ and to compute its residue at $x=1$. To this end, let us first consider the case when $\Delta \subset S \times S$ is cut out by a regular section $\sigma$ of a rank 2 vector bundle $\CV^\vee$ on $S \times S$. We have the Koszul complex of $\sigma$
\begin{equation}
\label{eqn:ses diag}
0\rightarrow \wedge^2\CV \rightarrow \CV \rightarrow \CO_{S \times S} \rightarrow \CO_{\Delta} \rightarrow 0
\end{equation}
If we write $[\CV] = L_1+L_2$ in $K_{S \times S}$, then we have
\begin{equation}
\label{eqn:euler ell}
\Delta_*(1) = \vartheta(L_1)\vartheta(L_2)
\end{equation}
and
$$
\zeta^\ellcoh(x) = \frac {\vartheta(xL_1)\vartheta(xL_2)}{\vartheta(x)\vartheta(xL_1L_2)}
$$
There is a general identity (where $L_1$ and $L_2$ are formal symbols)
\begin{equation}
\label{eqn:general identity}
\frac {\vartheta(xL_1)\vartheta(xL_2)}{\vartheta(x)\vartheta(xL_1L_2)}  = 1 + \vartheta(L_1)\vartheta(L_2) \cdot \frac {\text{quasi-periodic holomorphic function of }x}{\vartheta(x)\vartheta(xL_1L_2)}
\end{equation}
Because $\Delta$ is cut out by a regular section of $\CV^\vee$, we have
$$
\CV^\vee|_{\Delta} = \text{Nor}_{\Delta / S \times S}^\vee \cong \text{Tan}_S
$$
(where $\text{Nor}$ denotes the normal bundle) which implies that 
$$
L_1L_2|_\Delta = \det \CV |_\Delta = \det \text{Tan}_S^\vee = q
$$
Formulas \eqref{eqn:euler ell} and \eqref{eqn:general identity} imply \eqref{eqn:zeta want}, while \eqref{eqn:res zeta} follows from
$$
\underset{x=1}{\text{Res}} \ \zeta^{\ellcoh}(x) = - \frac {\vartheta(L_1)\vartheta(L_2)}{\vartheta(L_1L_2)} = \Delta_* \left(\frac {-1}{\vartheta(q)} \right)
$$
Outside of the special case when $\Delta$ is cut out by a regular section of a rank 2 vector bundle on $S \times S$, we may use deformation to the normal bundle to reduce to this special case (see \cite[Proposition 5.24]{N shuffle surf} for a version of this argument in $K$-theory).

\end{proof}

\medskip

\subsection{}
\label{sub:elliptic proof}

In what follows, we will use the obvious analogues of formulas \eqref{eqn:composition xy} and \eqref{eqn:composition yx} when discussing compositions of operators in elliptic cohomology. Let
\begin{equation}
\label{eqn:universal subset}
\ellcoh_{\CM \times S^k}^\circ \subseteq \ellcoh_{\CM \times S^k}
\end{equation}
denote the subset of all universal classes \eqref{eqn:universal class} (more precisely, universal classes form subsets not of $\ellcoh_{\CM \times S^k}$, but of various locally free rank 1 modules over $\ellcoh_{\CM \times S^k}$). We expect \eqref{eqn:universal subset} to be an equality; the standard way to prove such a statement is to use a resolution of the diagonal in $\CM \times \CM$ by universal classes.

\medskip

\begin{theorem}
\label{thm:elliptic}

We have the following equalities of operators $\eellcoh^\circ_{\CM} \rightarrow \eellcoh^\circ_{\CM \times \textcolor{red}{S} \times \textcolor{blue}{S}}$
\begin{equation}
\label{eqn:rel 1 elliptic}
\textcolor{red}{e(z)} \textcolor{blue}{e(w)} \tzeta^{\eellcoh}_{-}\left(\frac wz\right) =  \textcolor{blue}{e(w)}\textcolor{red}{e(z)}  \tzeta^{\eellcoh}_{+}\left(\frac zw\right) 
\end{equation}
\begin{equation}
\label{eqn:rel 2 elliptic}
 \textcolor{blue}{f(w)} \textcolor{red}{f(z)} \tzeta^{\eellcoh}_{-}\left(\frac wz\right) = \textcolor{red}{f(z)} \textcolor{blue}{f(w)} \tzeta^{\eellcoh}_{+}\left(\frac zw\right) 
\end{equation}
\begin{equation}
\label{eqn:rel 3 elliptic}
\textcolor{red}{h^\pm(z)} \textcolor{blue}{e(w)} =  \textcolor{blue}{e(w)} \textcolor{red}{h^\pm(z)} \frac {\zeta^{\eellcoh}\left(\frac zw\right)}{\zeta^{\eellcoh}\left(\frac wz\right)}
\end{equation}
\begin{equation}
\label{eqn:rel 4 elliptic}
\textcolor{blue}{f(w)} \textcolor{red}{h^\pm(z)} =\textcolor{red}{h^\pm(z)} \textcolor{blue}{f(w)} \frac {\zeta^{\eellcoh}\left(\frac zw\right)}{\zeta^{\eellcoh}\left(\frac wz\right)}
\end{equation}
\begin{equation}
\label{eqn:rel 5 elliptic}
\left [\textcolor{red}{f(z)}, \textcolor{blue}{e(w)} \right] = \delta \left(\frac zw\right)  \cdot \Delta_* \left( \frac{h^+(z) - h^-(w)}{\vartheta(q)}\right)
\end{equation}
as well as $[\textcolor{red}{h^\pm(z)}, \textcolor{blue}{h^{\pm'}(w)}] = 0$, for all $\pm,\pm' \in \{+,-\}$. 

\end{theorem}

\medskip

\noindent Our proof of formulas \eqref{eqn:rel 3 elliptic}--\eqref{eqn:rel 5 elliptic} establishes them as equalities of operators
$$
\ellcoh_{\CM} \rightarrow \ellcoh_{\CM \times \textcolor{red}{S} \times \textcolor{blue}{S}}
$$
Meanwhile, although our proof of formulas \eqref{eqn:rel 1 elliptic}--\eqref{eqn:rel 2 elliptic} only holds on the subsets of universal classes $\ellcoh^\circ \subseteq \ellcoh$, we believe that these relations actually hold on the whole elliptic cohomology rings.

\medskip

\begin{proof} We start with \eqref{eqn:rel 1 elliptic}, and leave the analogous formula \eqref{eqn:rel 2 elliptic} as an exercise to the interested reader. Let us apply formula \eqref{eqn:formula e universal} for any universal class $\Psi(\CU)$:
$$
\textcolor{blue}{e(w)} \cdot \Psi(\CU) =  \Psi(\CU + \textcolor{blue}{w} \CO_{\textcolor{blue}{\Delta}}) \vartheta \left(\frac {\textcolor{blue}{wq}}{\textcolor{blue}{\CU}} \right) \Big|_{1-p} 
$$
The right-hand side of the expression above lies in a certain locally free rank 1 module over $\ellcoh_{\CM \times \textcolor{blue}{S}}$, and the universal sheaf $\textcolor{blue}{\CU}$ (as well as the canonical class $\textcolor{blue}{q}$) corresponds to the latter copy of $\textcolor{blue}{S}$. Meanwhile, inside the universal class $\Psi(\CU)$, we write $\CU$ in black because it represents the universal sheaf on auxiliary copies of $S$, as in \eqref{eqn:universal series}. Applying \eqref{eqn:formula e universal} again to the identity above yields
\begin{equation}
\label{eqn:two e's}
\textcolor{red}{e(z)} \textcolor{blue}{e(w)} \cdot \Psi(\CU) = \vartheta \left( \frac {wq}{z} \CO^\vee_\Delta \right)
\end{equation}
$$
 \Psi(\CU + \textcolor{red}{z} \CO_{\textcolor{red}{\Delta}} + \textcolor{blue}{w} \CO_{\textcolor{blue}{\Delta}})  \vartheta \left(\frac {\textcolor{red}{zq}}{\textcolor{red}{\CU}} \right) \vartheta \left(\frac {\textcolor{blue}{wq}}{\textcolor{blue}{\CU}} \right) \Big|^{\textcolor{blue}{w}}_{1-p} \Big|^{\textcolor{red}{z}}_{1-p} 
$$
where the black copy of $\Delta$ above represents the diagonal in $\textcolor{red}{S} \times \textcolor{blue}{S}$, and the symbol $|_{1-p}^{\textcolor{blue}{w}} |_{1-p}^{\textcolor{red}{z}}$ means that we first expand in $w$ and then in $z$. However, note that
\begin{equation}
\label{eqn:dual identity}
q [\CO_\Delta^\vee] = [\CO_\Delta]  \in K_{\textcolor{red}{S} \times \textcolor{blue}{S}}
\end{equation}
implies $\vartheta \left( \frac {wq}{z} \CO^\vee_\Delta \right) =  \vartheta \left( \frac {w}{z} \CO_\Delta \right) = \zeta^{\ellcoh} \left(\frac wz\right)^{-1}$. Thus, formula \eqref{eqn:two e's} reads
$$
\textcolor{red}{e(z)} \textcolor{blue}{e(w)} \tzeta^{\ellcoh}_{-}\left(\frac wz\right) \cdot \Psi(\CU) = \vartheta \left(\frac {\textcolor{red}{zq}}{\textcolor{blue}{w}} \right) \vartheta \left(\frac {\textcolor{blue}{wq}}{\textcolor{red}{z}} \right)
$$
$$
 \Psi(\CU + \textcolor{red}{z} \CO_{\textcolor{red}{\Delta}} + \textcolor{blue}{w} \CO_{\textcolor{blue}{\Delta}}) \vartheta \left(\frac {\textcolor{red}{zq}}{\textcolor{red}{\CU}} \right) \vartheta \left(\frac {\textcolor{blue}{wq}}{\textcolor{blue}{\CU}} \right) \Big|^{\textcolor{blue}{w}}_{1-p} \Big|^{\textcolor{red}{z}}_{1-p} 
$$
Formula \eqref{eqn:rel 1 elliptic} follows from the fact that the right-hand side of the expression above is symmetric in $(\textcolor{red}{z},\textcolor{red}{S}) \leftrightarrow (\textcolor{blue}{w},\textcolor{blue}{S})$ (moreover, because there are no poles involving both $z$ and $w$, we can freely switch the order of the expansions $|_{1-p}^{\textcolor{red}{z}}$ and $|_{1-p}^{\textcolor{blue}{w}}$).

\medskip

\noindent We will now prove \eqref{eqn:rel 3 elliptic}, and leave the analogous formula \eqref{eqn:rel 4 elliptic} as an exercise to the reader. Recall that $\textcolor{blue}{e(w)}$ is given by the correspondence $\delta\left(\frac {L}w\right)$ on the locus
$$
\fZ = \Big\{ (\CF' \subset_{\textcolor{blue}{y}} \CF) \Big\}
$$
which is endowed with the line bundle $\CL$ (where $L = [\CL]$) and the map $p_S : \fZ \rightarrow \textcolor{blue}{S}$ that records the point $\textcolor{blue}{y}$. We have the short exact sequence
\begin{equation}
\label{eqn:ses colored}
0 \rightarrow \textcolor{red}{\CU'} \rightarrow \textcolor{red}{\CU} \rightarrow \CL \otimes \CO_{\Delta} \rightarrow 0
\end{equation}
on $\fZ \times \textcolor{red}{S}$, where $\textcolor{red}{\CU}$ and $\textcolor{red}{\CU'}$ are the universal sheaves keeping track of $\CF$ and $\CF'$ (respectively) as sheaves on $\textcolor{red}{S}$. In \eqref{eqn:ses colored}, we write $\CO_\Delta$ for the pull-back of the structure sheaf of the diagonal $\Delta \subset \textcolor{red}{S} \times \textcolor{blue}{S}$ under the map
$$
\fZ \times \textcolor{red}{S} \rightarrow \textcolor{red}{S} \times \textcolor{blue}{S}, \qquad (\CF' \subset_{\textcolor{blue}{y}} \CF) \times \textcolor{red}{x} \mapsto (\textcolor{red}{x}, \textcolor{blue}{y})
$$
With this in mind, note that the left-hand side of \eqref{eqn:rel 3 elliptic} is given by the following correspondence on $\fZ \times \textcolor{red}{S}$ 
\begin{align*}
\delta\left(\frac {L}{\textcolor{blue}{w}}\right) \frac {\displaystyle \vartheta \left( \frac {\textcolor{red}{zq}}{\textcolor{red}{\CU'}} \right)}{\displaystyle \vartheta \left( \frac {\textcolor{red}{z}}{\textcolor{red}{\CU'}} \right)} &= \delta\left(\frac {L}{\textcolor{blue}{w}}\right) \frac {\displaystyle \vartheta \left( \frac {\textcolor{red}{zq}}{\textcolor{red}{\CU}} \right)}{\displaystyle \vartheta \left( \frac {\textcolor{red}{z}}{\textcolor{red}{\CU}} \right)} \frac {\vartheta\left(\displaystyle -\frac {\textcolor{red}{zq}}L \CO^\vee_\Delta\right)}{\vartheta\left(\displaystyle  -\frac {\textcolor{red}{z}}L \CO^\vee_\Delta\right)} \\ &= \delta\left(\frac {L}{\textcolor{blue}{w}}\right) \frac {\displaystyle \vartheta \left( \frac {\textcolor{red}{zq}}{\textcolor{red}{\CU}} \right)}{\displaystyle \vartheta \left( \frac {\textcolor{red}{z}}{\textcolor{red}{\CU}} \right)} \frac {\vartheta\left(\displaystyle -\frac {\textcolor{red}{zq}}{\textcolor{blue}{w}} \CO^\vee_\Delta\right)}{\vartheta\left(\displaystyle  -\frac {\textcolor{red}{z}}{\textcolor{blue}{w}} \CO^\vee_\Delta\right)} \\
&= \delta\left(\frac {L}{\textcolor{blue}{w}}\right) \frac {\displaystyle \vartheta \left( \frac {\textcolor{red}{zq}}{\textcolor{red}{\CU}} \right)}{\displaystyle \vartheta \left( \frac {\textcolor{red}{z}}{\textcolor{red}{\CU}} \right)} \frac {\vartheta\left(\displaystyle -\frac {\textcolor{red}{z}}{\textcolor{blue}{w}} \CO_\Delta\right)}{\vartheta\left(\displaystyle  -\frac {\textcolor{blue}{w}}{\textcolor{red}{z}} \CO_\Delta\right)}
\end{align*}
(in the last equality, we used \eqref{eqn:dual identity} in the numerator, while in the denominator we used \eqref{eqn:change} together with the fact that $\CO_\Delta$ has rank 0 and determinant 1). The expression above is the correspondence that gives the right-hand side of \eqref{eqn:rel 3 elliptic}.

\medskip

\noindent Let us now turn to proving relation \eqref{eqn:rel 5 elliptic}. In Subsection \ref{sub:diagonal 1}, we showed that the commutator $[\textcolor{red}{f(z)}, \textcolor{blue}{e(w)}]$ is given by a correspondence supported on the locus $\{\textcolor{red}{x} = \textcolor{blue}{y}\}$ of $\CM \times \CM \times \textcolor{red}{S} \times \textcolor{blue}{S}$. However, a slight refinement of this argument (see the proof of Proposition 3.6 in \cite{N shuffle surf}) actually shows that the correspondence in question is supported on the smaller locus $\{\CF = \CF'\} \cap \{\textcolor{red}{x} = \textcolor{blue}{y}\}$. Therefore, we have
\begin{equation}
\label{eqn:comm}
[\textcolor{red}{f(z)}, \textcolor{blue}{e(w)}] = \Delta_* \Big( \text{multiplication by } \gamma \Big)
\end{equation}
for some elliptic cohomology class $\gamma$ on $\CM \times S$. It therefore remains to show that
\begin{equation}
\label{eqn:remains 0}
\gamma = \delta \left(\frac zw \right) \cdot  \frac 1{\vartheta(q)} \cdot \frac{\vartheta \left(\frac {zq}{\CU} \right)}{\vartheta \left(\frac {z}{\CU} \right)} \Big|_{1-p}
\end{equation}
Since $\Delta_*$ is injective (as it has a left-inverse given by projection onto the first factor $S \times S \rightarrow S$), formula \eqref{eqn:remains 0} is equivalent to
\begin{equation}
\label{eqn:remains 1}
\Delta_*(\gamma) = \delta \left(\frac zw \right) \cdot \Delta_* \left( \frac 1{\vartheta(q)} \cdot \frac{\vartheta \left(\frac {zq}{\CU} \right)}{\vartheta \left(\frac {z}{\CU} \right)} \Big|_{1-p} \right)
\end{equation}
In order to obtain $\Delta_*(\gamma)$, we simply apply equality \eqref{eqn:comm} to the unit class
$$
[\textcolor{red}{f(z)}, \textcolor{blue}{e(w)}]  \cdot 1 = \Delta_*(\gamma)
$$
and so it remains to show that
\begin{equation}
\label{eqn:remains 2}
[\textcolor{red}{f(z)}, \textcolor{blue}{e(w)}] \cdot 1 = \delta \left(\frac zw \right) \cdot \Delta_* \left( \frac 1{\vartheta(q)} \cdot \frac{\vartheta \left(\frac {zq}{\CU} \right)}{\vartheta \left(\frac {z}{\CU} \right)} \Big|_{1-p} \right)
\end{equation}
To prove \eqref{eqn:remains 2}, note that Proposition \ref{prop:e and f universal} gives us
\begin{align*}
&\textcolor{red}{f(z)} \cdot 1 =  \vartheta \left(- \frac {\textcolor{red}{z}}{\textcolor{red}{\CU}} \right) \Big|^{\textcolor{red}{z}}_{1-p} \\ 
&\textcolor{blue}{e(w)} \cdot 1 =  \vartheta \left(\frac {\textcolor{blue}{wq}}{\textcolor{blue}{\CU}} \right) \Big|^{\textcolor{blue}{w}}_{1-p}  
\end{align*}
(the color of the universal sheaf $\CU$, as well as the class $q = [\CK_S]$, denotes which copy of $\textcolor{red}{S} \times \textcolor{blue}{S}$ it corresponds to). Another application of Proposition \ref{prop:e and f universal} yields
\begin{align}
&\textcolor{red}{f(z)} \textcolor{blue}{e(w)} \cdot 1 = \frac {\displaystyle \vartheta \left(\frac {\textcolor{blue}{wq}}{\textcolor{blue}{\CU}} \right)}{\displaystyle \vartheta \left(\frac {\textcolor{red}{z}}{\textcolor{red}{\CU}} \right)} \vartheta \left( - \frac {wq}z \CO_\Delta^\vee\right) \Big|^{\textcolor{blue}{w}}_{1-p} \Big|^{\textcolor{red}{z}}_{1-p} \label{eqn:comp 1} \\
&\textcolor{blue}{e(w)} \textcolor{red}{f(z)} \cdot 1 = \frac {\displaystyle \vartheta \left(\frac {\textcolor{blue}{wq}}{\textcolor{blue}{\CU}} \right)}{\displaystyle \vartheta \left(\frac {\textcolor{red}{z}}{\textcolor{red}{\CU}} \right)} \vartheta \left( - \frac {w}{z} \CO_\Delta \right) \Big|^{\textcolor{red}{z}}_{1-p} \Big|^{\textcolor{blue}{w}}_{1-p} \label{eqn:comp 2} 
\end{align}
Because of formula \eqref{eqn:dual identity}, the right-hand sides of \eqref{eqn:comp 1} and \eqref{eqn:comp 2} only differ in the order in which we expand the variables $z$ and $w$, and so the difference $[\textcolor{red}{f(z)}, \textcolor{blue}{e(w)}]\cdot 1$ only has contributions from the residues of
$$
\vartheta \left( - \frac {wq}z \CO_\Delta^\vee\right) \stackrel{\eqref{eqn:dual identity}}=\vartheta \left( - \frac {w}z \CO_\Delta\right) = \zeta^{\ellcoh}\left(\frac wz\right)
$$
as $z$ passes over $w$. According to Lemma \ref{lem:integral}, there are only two such residues:

\medskip

\begin{itemize}[leftmargin=*]

\item The residue at $z = w$ yields 
$$
\delta \left(\frac z{w} \right) \cdot \Delta_* \left(\frac 1{\vartheta(q)} \cdot \frac{\vartheta \left(\frac {zq}{\CU} \right)}{\vartheta \left(\frac {z}{\CU} \right)} \Big|_{1-p}\right)
$$
which is equal to the right-hand side of \eqref{eqn:remains 2}. 

\medskip

\item The residue at $z = wq$ yields
$$
\delta \left(\frac z{wq} \right) \cdot \Delta_* \left( \text{const} \cdot \frac{\vartheta \left(\frac {z}{\CU} \right)}{\vartheta \left(\frac {z}{\CU} \right)} \Big|_{1-p}\right)
$$
where $\text{const}$ denotes minus the residue of $\zeta^{\ellcoh}(x)$ at $x=q^{-1}$. The expression above is 0 because the two $\vartheta$ in the right-hand side cancel out, and $\text{const}|_{1-p} = 0$.

\end{itemize}

\end{proof}

\medskip

\subsection{}
\label{sub:ko}

We will now connect the formulas in Theorem \ref{thm:elliptic} with the elliptic quantum toroidal algebra studied in \cite{ko}. The main modification we make from the elliptic curve context in the previous Subsections is that we will think of $p$ as being infinitesimally small (instead of being a fixed complex number) and so we will expand all of our formulas as power series in $p$. Consider the meromorphic functions
\begin{equation}
\label{eqn:def zeta ell}
\zeta^{E}(x) = \frac {\vartheta(xq_1)\vartheta(xq_2)}{\vartheta(x)\vartheta(xq)}
\end{equation}
\begin{equation}
\label{eqn:def zeta ell tilde}
\tzeta^{E}(x) = \zeta^{E}(x)\vartheta(xq)\vartheta(x^{-1}q) = \frac {\vartheta(xq_1)\vartheta(xq_2)\vartheta(x^{-1}q)}{\vartheta(x)}
\end{equation} 
where $q_1,q_2$ are formal symbols, and $q = q_1q_2$.

\medskip

\begin{definition}
\label{def:ell}

(\cite{ko}) Elliptic quantum toroidal $\fgl_1$ is the algebra
$$
U_{q_1,q_2,p}(\ddot{\fgl}_1) = \BZ[q^{\pm 1}_1,q^{\pm 1}_2][[p]] \Big\langle e_n, f_n, h_n^\pm \Big \rangle_{n \in \BZ} \Big / \text{relations \eqref{eqn:rel 1 ell}--\eqref{eqn:rel 5 ell}}
$$
The defining relations are best written in terms of the generating series
$$
e(z) = \sum_{n=-\infty}^{\infty} \frac {e_n}{z^n}, \qquad f(z) = \sum_{n=-\infty}^{\infty} \frac {f_n}{z^n}, \qquad h^\pm(z) = \sum_{n=-\infty}^{\infty} \frac {h^\pm_n}{z^n}
$$
and take the form
\begin{equation}
\label{eqn:rel 1 ell}
e(z)e(w) \tzeta^{E}\left(\frac wz \right) = e(w)e(z) \tzeta^{E}\left(\frac zw \right)
\end{equation}
\begin{equation}
\label{eqn:rel 2 ell}
f(w) f(z) \tzeta^{E}\left(\frac wz \right) = f(z)f(w) \tzeta^{E}\left(\frac zw \right)
\end{equation}
\begin{equation}
\label{eqn:rel 3 ell}
h^\pm(z)e(w) =  e(w)h^\pm(z) \frac {\zeta^{E}\left(\frac zw \right)}{\zeta^{E}\left(\frac wz \right)}
\end{equation}
\begin{equation}
\label{eqn:rel 4 ell}
f(w) h^\pm(z) = h^\pm(z) f(w) \frac {\zeta^{E}\left(\frac zw \right)}{\zeta^{E}\left(\frac wz \right)}
\end{equation}
\begin{equation}
\label{eqn:rel 5 ell}
\left [f(z), e(w) \right] = \frac {\vartheta(q_1)\vartheta(q_2)}{\vartheta(q)} \cdot \delta \left(\frac zw\right) \Big( h^+(z) - h^-(w) \Big)
\end{equation}
(where $\delta(x) = \sum_{n=-\infty}^{\infty} x^n$), as well as $[h^\pm(z),h^{\pm'}(w)] = 0$ for all $\pm,\pm' \in \{+,-\}$.

\end{definition}

\medskip

\noindent Note that beside the generators and relations above, the elliptic quantum toroidal algebra of \cite{ko} has one more central element $\gamma$ (which is always set equal to 1 in geometric representations) and one more cubic relation involving the $e_n$'s and $f_n$'s (respectively). We do not include the latter relation in order to streamline the presentation; it is straightforward to check that the cubic relations between the operators \eqref{eqn:e elliptic series} and \eqref{eqn:f elliptic series} (respectively) hold by using formulas \eqref{eqn:formula e universal} and \eqref{eqn:formula f universal}. 

\medskip 

\noindent Comparing \eqref{eqn:rel 1 elliptic}--\eqref{eqn:rel 5 elliptic} with Definition \ref{def:ell} allows us to phrase Theorem \ref{thm:elliptic} as yielding an action (in the sense analogous to Definition \ref{def:action})
\begin{equation}
\label{eqn:ell coh acts}
U_{q_1,q_2,p}(\ddot{\fgl}_1) \curvearrowright \ellcoh^\circ_{\CM}
\end{equation}
We do not make this into a precise Theorem for rather pedantic reasons, such as the fact that \cite{ko} define the elliptic quantum group over power series in $p$, while $\ellcoh_{\CM}$ is defined for fixed $p$ (this discrepancy is important when discussing the integrality of the algebras involved). However, \eqref{eqn:ell coh acts} is the motivation for the present paper.

\medskip

\begin{remark}

In \cite{yz}, the authors define an elliptic quantum group in a scheme-theoretic way (in the related setup where moduli spaces of sheaves on surfaces are replaced by Nakajima quiver varieties), following the axiomatic treatment of \cite{gkv} and the theory of cohomological Hall algebras. To obtain explicit algebras out of their general construction, one needs to choose a ``basis" of meromorphic functions on the elliptic curve $E$. The choice we made in the present paper is that of the monomials $\{z^n\}_{n \in \BZ}$, which led to the particular algebra of Definition \ref{def:ell}. However, other choices (such as
$$
\left\{\frac {\partial^n}{\partial z^n} \left( \frac {\vartheta(zu)}{\vartheta(z)\vartheta(u)} \right) \right\}_{n \geq 0, u \in \BC^*}
$$
that was considered in \cite{yz}) are equally valid, and they produce different generators-and-relations realizations of elliptic quantum groups.

\end{remark}

\end{document}